\documentclass[12pt]{article}
\usepackage[a4paper,left=2cm,right=1.2cm,top=3cm,bottom=4cm]
{geometry}
\usepackage{amsmath,amstext, amsthm, amsbsy,amssymb,marvosym,fancyhdr,graphicx,amscd,amsfonts,latexsym,delarray,stackrel,lineno,color,cite,appendix}
\usepackage{wasysym}
\usepackage{mathrsfs}
\usepackage{caption}
\usepackage{subcaption}
\usepackage{srcltx}
\usepackage{mathrsfs}
\usepackage{float} 
\usepackage[all]{xy}
\usepackage{t1enc}
\usepackage{mathrsfs}
\usepackage{pifont}

\usepackage{mathpple}
\usepackage[T1]{fontenc}

\newcommand*\patchAmsMathEnvironmentForLineno[1]{%
  \expandafter\let\csname old#1\expandafter\endcsname\csname #1\endcsname
  \expandafter\let\csname oldend#1\expandafter\endcsname\csname end#1\endcsname
  \renewenvironment{#1}%
     {\linenomath\csname old#1\endcsname}%
     {\csname oldend#1\endcsname\endlinenomath}}%
\newcommand*\patchBothAmsMathEnvironmentsForLineno[1]{%
  \patchAmsMathEnvironmentForLineno{#1}%
  \patchAmsMathEnvironmentForLineno{#1*}}%
\AtBeginDocument{%
\patchBothAmsMathEnvironmentsForLineno{equation}%
\patchBothAmsMathEnvironmentsForLineno{align}%
\patchBothAmsMathEnvironmentsForLineno{flalign}%
\patchBothAmsMathEnvironmentsForLineno{alignat}%
\patchBothAmsMathEnvironmentsForLineno{gather}%
\patchBothAmsMathEnvironmentsForLineno{multline}%
}
\definecolor{Green}{rgb}{0,1,0}
\definecolor{Blue}{RGB}{0,0,191}
\definecolor{mathmodecolor}{RGB}{0,102,0}
\definecolor{keywordcolor}{RGB}{0,51,151}
\definecolor{sourcebackgroundcolor}{RGB}{255,247,223}
\definecolor{unixagred}{RGB}{255,0,0}
\definecolor{lightgray}{RGB}{191,191,191}
\definecolor{green}{RGB}{1,191,191}

\newtheorem{thm}{Theorem}[section]
\newtheorem*{thm*}{Theorem}
\newtheorem{prop}[thm]{Proposition}

\newtheorem{lem}[thm]{Lemma}

\newtheorem{defn}[thm]{Definition}
\newtheorem{rem}[thm]{Remark}

\def\Aut{{\rm Aut}}

\def\Dom{{\rm Dom}}
\def\End{{\rm End}}

\def\Hom{{\rm Hom}}
\def\id{{\rm id}}

\def\Ker{{\rm Ker}}

\def\Spec{{\rm Spec\,}}
\def\Sp{{\rm Spec}\,}

\def\A{{\mathbb A}}
\def\B{{\mathbb B}}
\def\C{{\mathbb C}}
\def\F{{\mathbb F}}
\def\K{{\mathbb K}}
\def\N{{\mathbb N}}
\def\P{{\mathbb P}}
\def\Q{{\mathbb Q}}
\def\R{{\mathbb R}}
\def\Z{{\mathbb Z}}

\newcommand{\cff}{{\it cf.}~}

\def\aarith{{\mathfrak A}}

\def\ord{{\rm Ord}}
\def\ntp{\N_{(p)}^{\times}}

\def\Codom{{\mbox{Codom}}}
\def\cdiv{{\mathcal{C}a\mathcal{C}\ell}}

\def\cA{{\mathcal A}}

\def\cE{{\mathcal E}}

\def\cF{{\mathcal F}}

\def\cK{{\mathcal K}}

\def\cO{{\mathcal O}}
\def\cP{{\mathcal P}}

\def\cR{{\mathcal R}}

\def\cU{{\mathcal U}}
\def\cV{{\mathcal V}}

\def\cZ{{\mathcal Z}}
\def\scal{{(\rnt,\cO)}}
\def\scal1{{\hat \aarith}}

\def\qqq{\,,\,~\forall}

\newcommand{\ie}{{\it i.e.\/}\ }
\newcommand{\eg}{{\it e.g.\/}\ }
\newcommand{\cf}{{\it cf.}}

\def\supp{{\rm Support }}

\def\id{{\mbox{Id}}}

\def\dim{{\mbox{dim}}}
\def\cdim{{{\mbox{Dim}_\R}}}
\def\tdim{{{\mbox{dim}_{\rm top}}}}

\def\Hom {{\mbox{Hom}}}

\def\End{{\mbox{End}}}

\def\ord{{\rm Ord}}

\def\ffp{\mathfrak{p}}
\def\ffq{\mathfrak{q}}

\def\mc{multiplicatively cancellative }

\def\rma{\R_{\rm max}}
\def\sss{{\mathbb S}}

\def\zmax{{\Z_{\rm max}}}

\def\rmax{\R_+^{\rm max}}

\def\fr{{\rm Fr}}

\def\Se{\frak{ Sets}}

\def\sh{\mathfrak{Sh}}

\newcommand{\nil}[1]{}

\def\nt{\N^{\times}}
\def\wnt{{\widehat{\N^{\times}}}}
\def\rnt{{[0,\infty)\rtimes{\N^{\times}}}}

\def\div{{\rm Div}}

\def\aarith{{\mathscr A}}
\def\scal{{(\rnt,\cO)}}
\def\scal1{{\hat \aarith}}
\def\scal2{{\mathscr S}}
\title
{Geometry of the scaling site}
\date{}

\author{Alain Connes, Caterina Consani\thanks{Partially supported by the Simons Foundation collaboration grant n.~353677. This author would also like to thank the Coll\`ege de France for some financial support.}}

\begin{document}

\maketitle

\begin{abstract}
In this paper we construct the scaling site $\scal2$ by implementing the extension of scalars on the arithmetic site $\aarith$, from the smallest Boolean semifield $\B$ to the tropical semifield $\rmax$. The obtained semiringed topos is the Grothendieck topos $\rnt$, semi-direct product of the Euclidean half-line and the monoid $\nt$ of positive integers acting by multiplication, endowed with the structure sheaf of piecewise affine, convex functions with integral slopes. We show that pointwise $\rnt$ coincides with the adele class space of $\Q$ and that this latter space inherits the  geometric structure of a tropical curve. We restrict  this construction to the periodic orbit of the scaling flow associated to each prime $p$  and obtain a quasi-tropical structure which turns this orbit into a variant $C_p=\R_+^*/p^\Z$ of the classical Jacobi description $\C^*/q^\Z$ of an elliptic curve. On $C_p$, we develop the theory of Cartier divisors, determine the structure of the quotient $\div(C_p)/\cP$ of the abelian group of divisors by the subgroup of principal divisors, develop the theory of theta functions, and prove the Riemann-Roch formula which involves real valued dimensions, as in the type II index theory. 
We show  that one would have been led to the same definition of $\scal2$  by analyzing the well known results on the localization of zeros of analytic functions involving Newton polygons in the non-archimedean case and the Jensen's formula in the complex case.
\end{abstract}
\vspace{0.1in}

\tableofcontents

\section{Introduction}
\begin{figure}
\begin{center}
\includegraphics[scale=0.3]{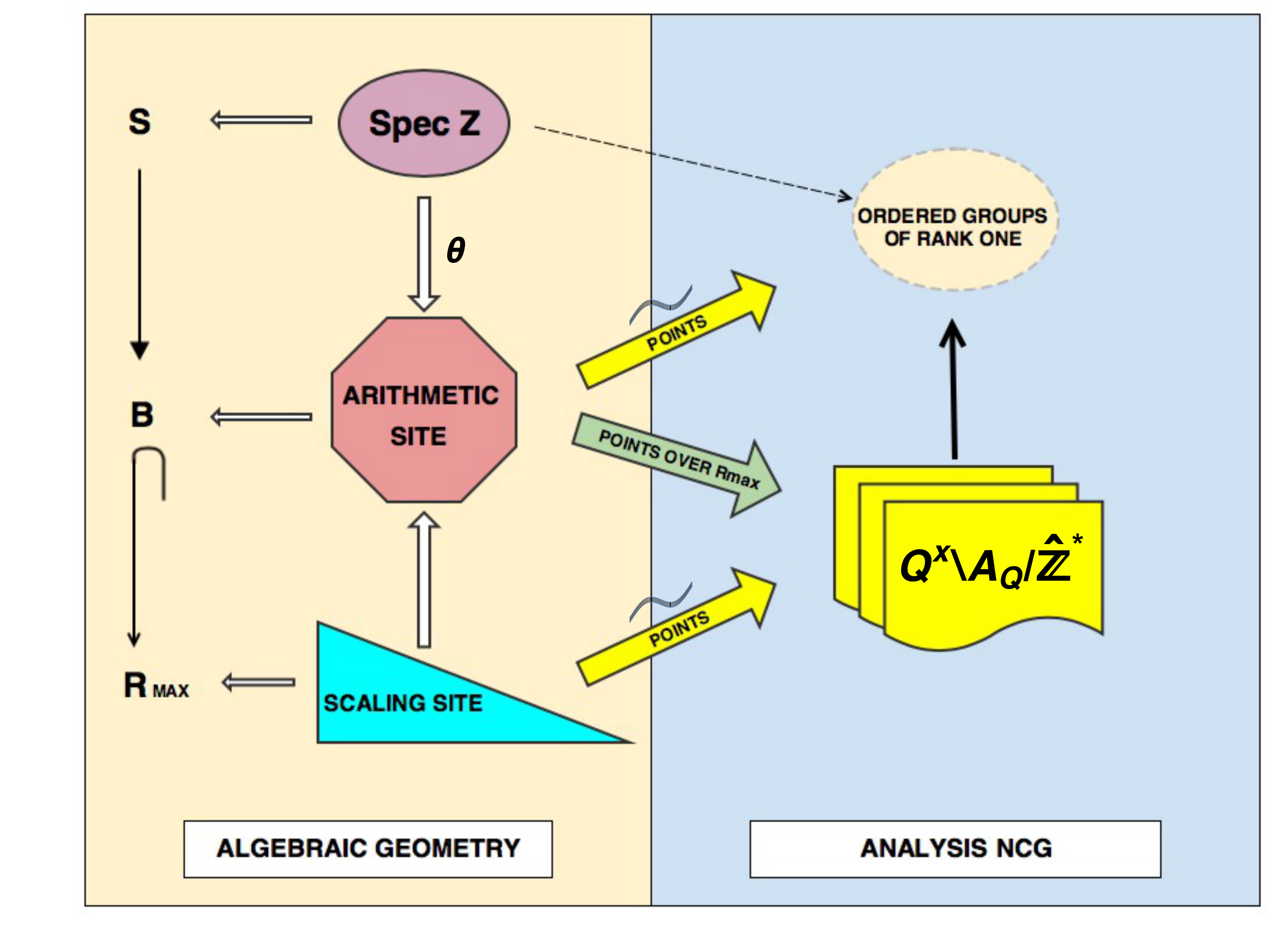}
\end{center}
\caption{The extension of scalars from $\aarith$ to $\scal2$ \label{scalingpic} }
\end{figure}

In this article we show that the adele class space\footnote{More specifically the sector $\Q^\times\backslash\A_\Q/\hat\Z^*$ corresponding to the trivial Grossencharacter.} of $\Q$  admits a natural structure of tropical curve. We follow the strategy outlined in \cite{CC,CC1} and investigate the algebraic geometric structure of the 
 Scaling Site\footnote{These results have been announced in \cite{CC2}.} $\scal2$   obtained from the arithmetic site $\aarith$  by extension of scalars from the Boolean semifield $\B$ to the tropical semifield $\rmax$ (\cf~Figure~\ref{scalingpic}). As a Grothendieck topos $\scal2$ is described as $\rnt$: the topos of $\nt$-equivariant sheaves (of sets) on the half-line  and our first result (Theorem \ref{scaltop}) states that the isomorphism classes of points  of this topos form the basic sector of the adele class space of $\Q$
\begin{thm*}\label{scaltopintro}The  space of  points of the topos $\rnt$
is canonically isomorphic to $\aarith(\rmax)\simeq \Q^\times\backslash\A_\Q/\hat\Z^*$.
\end{thm*}
This result provides the missing geometric structure on the adele class space since the topos $\rnt$ inherits, from its construction by extension of scalars, a natural sheaf $\cO$ of regular functions. We call  {\em Scaling Site}  the semi-ringed topos 
\[
\scal2:=\left(\rnt,\cO\right)
\]
 so obtained.  The sections of the sheaf $\cO$ are convex,  piecewise affine functions with integral slopes. In  Appendix \ref{apptropic} we review the well known results on the localization of zeros of analytic functions  showing  in which sense the tropical half-line $(0,\infty)$, endowed with the  sheaf of convex  piecewise affine functions with integral slopes, provides a suitable structure  for the localization (both in the archimedean and non-archimedean case) of zeros of analytic functions in the punctured  unit disk. The new component supplied with the scaling site is the action of $\nt$ by multiplication on the tropical half-line $[0,\infty)$. On  analytic functions this action is given by the transformation $f(z)\mapsto f(z^n)$ for $n\in \nt$, \ie  the action of the degree $n$ endomorphism $z\mapsto z^n$ on the punctured  unit disk. \newline 
 The structure sheaf $\cO$ of $\scal2$   is a sheaf of semirings 
 of ``characteristic one"  (\ie of semirings in which $1+1=1$) and the naturalness of this structure is justified at the conceptual level (see Appendix  \ref{appzmax}) by two  facts. First, the endomorphisms of any object admit a {\em canonical} structure of  semiring in any category  with finite products and coproducts when the canonical morphisms from coproducts to products are isomorphisms. Second, passing from rings to semirings only adds one more object to the list of finite fields $\F_q$, namely   the Boolean semifield $\B$, and only one object to the list of fields  whose multiplicative group is cyclic, \ie the semifield $\zmax$ whose multiplicative group is infinite cyclic. Both $\B$ and $\zmax$ are semirings of characteristic one.

 In section \ref{sectsheafstalks} we describe the stalks of the structure sheaf $\cO$  and show (Theorem \ref{structure3}) that the  points $\scal2(\rmax)$ of the scaling site defined over $\rmax$ coincide with the points 
 $\aarith(\rmax)$ of the arithmetic site defined over the same semifield. As stated in \cite{CC,CC1} a long term goal of this project is to develop  an adequate version of the Riemann-Roch theorem in characteristic $1$, suitable to transplant the pRH proof of Weil to the Riemann zeta function. In this paper we test this idea by restricting our geometric structure to the periodic orbits of the scaling flow, \ie to the points of $\rnt$ over the image of $\Spec\Z$ (\cf~Figure~\ref{scalingpic} and \cite{CC1}, \S 5.1). We find that for each prime $p$ the corresponding circle of length $\log p$ is endowed with a quasi-tropical structure which turns this orbit into a variant $C_p=\R_+^*/p^\Z$ of the classical Jacobi description $\C^*/q^\Z$ of an elliptic curve. The structure sheaf $\cO_p$ of $C_p$ is obtained by restriction of $\cO$ to $C_p$ and its sections are periodic functions $f(p\lambda)=f(\lambda)$, $\lambda \in \R_+^*$, which are convex, piecewise affine and whose derivatives take values in the group $H_p\subset \R$ of rational numbers with denominators a power of $p$. When suitably understood in conceptual terms using Cartier divisors, the notions of rational functions, divisors, etc. on $C_p$ are all meaningful. The global rational functions form a semifield $\cK(C_p)$ (of characteristic one). A new feature of this construction is that the degree of a divisor can  be any real number. We introduce an invariant $\chi(D)\in \Z/(p-1)\Z$ for divisors $D$ on $C_p$ and determine, in Theorem \ref{thmjaccp},  the precise structure of the quotient   $\div(C_p)/\cP$ of the abelian group of divisors by the subgroup of principal divisors
\begin{thm*}\label{thmjaccpintro}  The map $(\deg,\chi)$ is an isomorphism of abelian groups 
$$
(\deg,\chi):\div(C_p)/\cP\to \R\times (\Z/(p-1)\Z).
$$
\end{thm*}

 We develop, in analogy with the non-archimedean version established in \cite{Tate}, the theory of theta functions on $C_p$, starting with 
  the following infinite sums as the analogues\footnote{We use the notation $x\vee y$ for the max of two real numbers} of the infinite products defining classical theta functions
 $$
\theta(\lambda):=\sum_0^\infty \left(0 \vee (1-p^{m}\lambda)\right)
+\sum_1^\infty \left(0 \vee (p^{-m}\lambda-1)\right).
$$
We define $\theta$-functions $\Theta_{h,\mu}$ for $h\in H_p$, $h>0$  and $\mu\in \R_+^*$. They are obtained by applying to the basic theta function $\theta(\lambda)$ defined above the symmetries associated to the various incarnations (arithmetic, relative, absolute, geometric) of the ``Frobenius" operator in this context. This part is discussed in details in \S \ref{sectsymm}.   The main output (\cf~Theorem \ref{thmtheta1}) is provided by the following 
\begin{thm*} Any function $f\in \cK(C_p)$ is canonically expressed in terms of theta 
functions associated to the principal divisor of $f$, and a constant $c\in\R$
$$
f(\lambda):=\sum_i \Theta_{h_i,\mu_i}(\lambda)-\sum_j \Theta_{h'_j,\mu'_j}(\lambda)-h\lambda +c.
$$
\end{thm*} 
 For each divisor $D$ on $C_p$ we define the corresponding Riemann-Roch problem with solution space $H^0(D):=H^0(C_p,\cO(D))$. We introduce the continuous dimension $\cdim(H^0(D))$ of this $\rma$-module using a limit of normalized topological dimensions and find that $\cdim(H^0(D))$ is a real number. The topological dimension used in this part is the Lebesgue covering dimension which assigns to any topological space $X$ an integer $\tdim(X)\in \{0,\ldots ,\infty\}$ counting the minimal overlap of refinements of open covers. The appearance of arbitrary positive real numbers as continuous dimensions of $\cdim(H^0(D))$ is due to the density in $\R$ of the subgroup $H_p\subset \Q$ and the fact that continuous dimensions are defined as limits
 $$
 \cdim(H^0(D)):=\lim_{n\to \infty} p^{-n}\tdim(H^0(D)^{p^n})
 $$
  of normalized dimensions $p^{-n}\tdim(H^0(D)^{p^n})$ where $H^0(D)^{p^n}$ is a natural filtration of $H^0(D)$ involving the $p$-adic norms of the derivatives. We interpret this result as the characteristic $1$ counterpart of the statement for matroid  $C^*$-algebras and the type II normalized traces as in \cite{dix}. The continuous dimensions which affect arbitrary positive {\em real} values appear when passing to the von Neumann algebra of type II obtained as the weak closure of the $C^*$-algebra using the trace to perform the completion. Finally, in Theorem \ref{RRperiodic} we prove  that the  Riemann-Roch formula holds for $C_p$
 \begin{thm*}\label{RRperiodicintro}
Let $D\in \div(C_p)$ be a divisor, then the limit 
$
\cdim(H^0(D))
$ exists and  
one has the Riemann-Roch formula: 
$$
\cdim(H^0(D))-\cdim(H^0(-D))=\deg(D)\qqq D\in \div(C_p).
$$	
\end{thm*}

By comparing the periodic orbit $C_p$ with a tropical elliptic curve and  our Riemann-Roch theorem with the tropical Riemann-Roch theorem of \cite{BN,GK,MZ} and its variants we find several fundamental differences. First, for an elliptic tropical curve   $C$ given by a circle of length $L$, the structure of the group $\div(C)/\cP$ of divisor classes is inserted into an exact sequence of the  form
$
0\to \R/L\Z\to \div(C)/\cP\stackrel{\deg}{\to} \Z\to 0
$ (\cf  \cite{MZ}\!\! ),
while for the periodic orbit $C_p$ the group of divisor classes is $\div(C_p)/\cP\simeq 
\R\times (\Z/(p-1)\Z)$. The second fundamental difference is the appearance of continuous dimensions in our Riemann-Roch theorem. 
The source for these differences is seen when one compares the structure sheaf of $C_p$
with that of the elliptic tropical curve $C:=\R/L\Z$, $L=\log p$. Let us use for $C_p$ the variable $u=\log\lambda$, so that  the periodicity condition $f(px)=f(x)$ becomes translation invariance by $\log p$. Then  the local sections of the structure sheaf of $C_p$ are in particular  piecewise affine in the parameter $\lambda$
and this condition  is expressed, in the variable $u$, by the piecewise vanishing of $\Delta_2f$, where $\Delta_2$ is the elliptic translation invariant operator 
\begin{equation}\label{maineq}
\Delta_2=\lambda^2\left(\frac{\partial}{\partial \lambda}\right)^2, \ \ \Delta_2(f):=(D_u^2-D_u)f, \  \  \ D_u:=\frac{\partial}{\partial u}.
\end{equation}
On the other hand, the local sections of the structure sheaf  of $C$ are in particular  piecewise affine in the parameter $u$, and this condition  is expressed, in the variable $u$, by the piecewise vanishing of $D_u^2\, f$. Thus one readily sees that the difference between the two sheaves is due to the presence of the sub-principal term $- D_u$ in \eqref{maineq}.

\subsection*{Notations} For any abelian ordered group $H$ we denote by $H_{\rm max}=H\cup \{-\infty\}$  the semifield obtained from $H$ by applying the max-plus construction \ie the addition is given by the max and the multiplication by the addition in $H$. Since $\R_{\rm max}$ is isomorphic to $\rmax$ by the exponential map (\cf \cite{Gaubert}) we shall pass freely from the ``additive" notation $\rma$ to the ``multiplicative" one $\rmax$. 
\vspace*{-.5cm}

\section{The topos $\rnt$}

In this section we define the topos underlying the scaling site $\scal2$  as a Grothendieck site, \ie as a small category $\mathscr C$ endowed with a Grothendieck topology $J$. In \S \ref{sectextS} we shortly explain  its structure as naturally arising  from the arithmetic site  $\aarith$ by extension of scalars from $\B$ to $\rmax$. In \S \ref{sectcatC} we provide the definition of the small category $\mathscr C$ and in \S \ref{sectGtop} we describe its Grothendieck topology.

\subsection{Extension of scalars}\label{sectextS}

The arithmetic site $\aarith$ of \cite{CC,CC1} is defined using the action of $\nt$ by Frobenius endomorphisms $\fr_n$ on the semifield $\zmax$ of characteristic one. To define the extension of scalars from $\B$ to $\rmax$ we consider the semiring $\cR(\Z)=\zmax\hat\otimes_\B\rma$
obtained as the \mc reduction of the tensor product $\zmax\otimes_\B\rma$ and we endow $\cR(\Z)$ with the $\rma$-linear endomorphisms $\fr_n\otimes \id$. Then, by applying  the Legendre transform we identify $\cR(\Z)$ with the semiring of  convex piecewise affine functions on 
$\R_+$ with slopes in $\Z\subset \R$ and only finitely many discontinuities of the derivative. These functions are endowed with the pointwise 
operations of functions taking values in $\rma$.

The operation of reduction from $\zmax\otimes_\B\rma$ to $\cR(\Z)=\zmax\hat\otimes_\B\rma$ is obtained as described in \cite{CC1} Lemma 6.20 and Proposition 6.21. More precisely, the 
 elements of $\cR(\Z)$ are given by the convex hull $C$ of the union of finitely many quadrants of the form $(x_j,y_j)-Q$, for $Q=\R_+\times \R_+$, where $x_j\in \Z$ and $y_j\in \R$. To determine this convex hull it is enough to know which half planes $P\subset \R^2$ contain it, and any such half-plane has the form 
$$
P_{\lambda,u}:=\{(x,y)\in\R^2\mid \lambda x+y\leq u\}, \qquad P^v:=\{(x,y)\in\R^2\mid x\leq v\}
$$
where $\lambda \in \R_+$ and $u,v\in\R$.
Thus $C$ is uniquely determined by the function 
\begin{equation*}\label{hdefn0}
\ell_C(\lambda):=\inf \{u\in\R\mid C\subset P_{\lambda,u}\}
\end{equation*}
and  this function is given in terms of the finitely many vertices $(x_j,y_j)$ of the  polygon $C$ by the formula 
\begin{equation}\label{hdefn}
\ell_C(\lambda)=\max_j \lambda x_j+y_j.
\end{equation}
\begin{figure}%[H]
%\centering
\begin{subfigure}{.5\textwidth}
  \centering
  \includegraphics[scale=0.3]{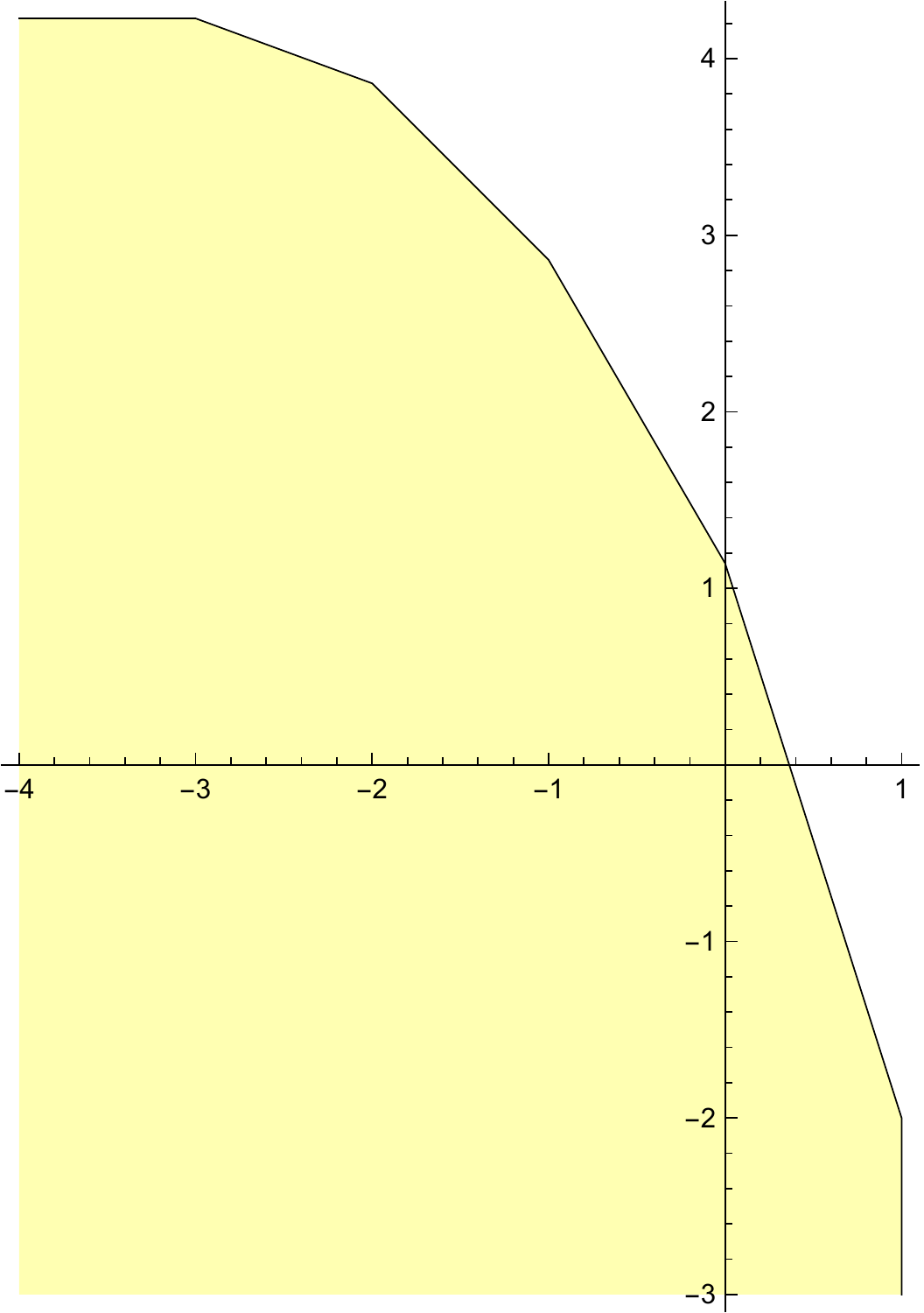}
  \caption{An element $C$ of $\cR(\Z)$}
  \label{elem}
\end{subfigure}
\begin{subfigure}{.5\textwidth}
  %\centering
  \includegraphics[scale=0.5]{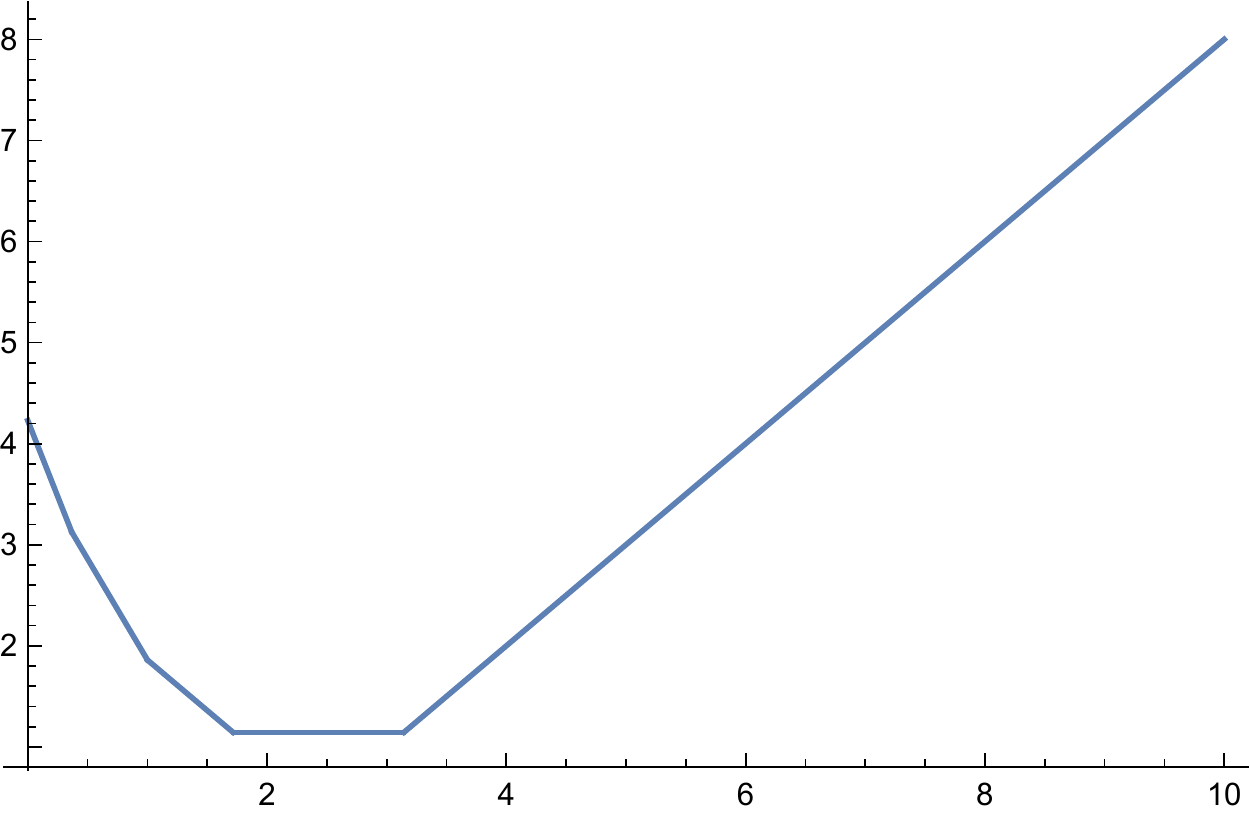}
  \caption{The Legendre transform $\ell_C(\lambda)$}
  \label{lege}
\end{subfigure}
\end{figure}

Note that $\ell_C(\lambda)$ is meaningful also for $\lambda=0$ and that  $\lim_{\lambda\to \infty}\, \ell_C(\lambda) /\lambda =\max x_j=\inf \{v\in \R\mid C\subset P^v\}$.
 One then obtains the required identification
\begin{lem} \label{legendrelem} 
The map $L:C\mapsto \ell_C$ where $\ell_C(\lambda)$, for $\forall \lambda\in \R_+$ has been defined in \eqref{hdefn}, is an isomorphism of $\zmax\hat\otimes_\B\rma$ with the semiring $\cR(\Z)$ of continuous\footnote{Continuity is automatic for convex functions on open intervals, see \cite{Rudin}, Theorem 3.2.} convex piecewise affine functions on 
$\R_+$ with slopes in $\Z\subset \R$ and only finitely many discontinuities of the derivative. These functions are endowed with the pointwise 
operations of functions with values in $\rma$.
\end{lem}
\proof It follows by the formula \eqref{hdefn} that the function $\ell_C$ belongs to $\cR(\Z)$ since the slopes $x_j\in \Z$, and the discontinuities of the derivative only occur when there exists a pair of vertices $(x_i,y_i)\neq (x_j,y_j)$ such that $\lambda(x_i-x_j)=y_j-y_i$. Moreover $\ell_C$ is convex  by construction as a supremum of finitely many affine functions. The map $C\mapsto \ell_C$ is evidently injective. The surjectivity is a consequence of the fact that  an element $f\in \cR(\Z)$ is a finite supremum of affine functions $\lambda\mapsto a\lambda+b$ with $a\in \Z$ and $b\in \R$, and thus of the form $\ell_C$ where $C$ is the convex hull of the union of quadrants $(a,b)-Q$. The  pointwise 
operations of functions with values in $\R_{\rm max}$ are given for the addition by the rule $(f,g)\mapsto f\vee g$, where $(f\vee g)(\lambda):=f(\lambda)\vee g(\lambda)=\max\{f(\lambda), g(\lambda)\}$. This corresponds to the $\max$ in \eqref{hdefn} and thus to the convex hull of the union in terms of the convex sets $C$. This shows that the map $L:C\mapsto \ell_C$ is additive. It is also multiplicative, \ie $\ell_{C+C'}=\ell_C+\ell_{C'}$. This follows using \eqref{hdefn} and the identity   $\max(A+B)=\max(A)+\max(B)$ holding for any two finite subsets $A,B\subset \R$.\endproof
\begin{prop}\label{propextscal} $(i)$~Under the isomorphism $L$ the endomorphism $\fr_n\otimes \id$ of $\zmax\hat\otimes_\B\rma$ corresponds to the action of $\nt$ on $\R_+$ by multiplication.\newline
$(ii)$~The following map identifies the half line $[0,\infty)$ with the space of characters of $\cR(\Z)$
$$
[0,\infty)\ni\lambda\mapsto \chi_\lambda \in \Hom_{\rma}(\cR(\Z),\rma), \qquad \chi_\lambda(f):=f(\lambda).
$$

\end{prop}
\proof $(i)$~We use multiplicative notations both for $\zmax$ and $\rma$ and represent elements of
$\zmax\hat\otimes_\B\R_{\rm max}$ in terms of finite sums  $\sum q^{x_j}\otimes_\B q^{y_j}$ with $x_j\in \Z$ and $y_j\in \R$. In these terms 
 the isomorphism $L:\zmax\hat\otimes_\B\R_{\rm max}\to \cR(\Z)$ is such that 
$$
L(\sum q^{x_j}\otimes_\B q^{y_j})(\lambda)=\max\{\lambda x_j+y_j\}
$$
With $X=\sum q^{x_j}\otimes_\B q^{y_j}$ one has 
$$
L((\fr_n\otimes \id)(X)(\lambda)=
L(\sum q^{nx_j}\otimes_\B q^{y_j})(\lambda)=\max\{n\lambda x_j+y_j\}=L(X)(n\lambda).
$$

$(ii)$~Let $\iota\in \cR(\Z)$ be the function $\iota(\lambda)=\lambda$.  An element $\rho\in \Hom_{\rma}(\cR(\Z),\rma)$ is uniquely specified  by  $\rho(\iota)\in \rma$ and since  $0\vee \iota=\iota$ and $\rho(0)=0$ (as the morphism $\rho$ preserves the multiplicative unit) one has  $\rho(\iota)=\lambda\in [0,\infty)\subset \rma$. By multiplicativity one gets $\rho(k\iota)=k\lambda$ for any $k\in \Z$ and by $\rma$-linearity that $\rho(k\iota +y)=k\lambda+y$ for $y\in \R$. By additivity one then gets that for any $f=\vee(x_j\iota+y_j)\in \cR(\Z)$, $\rho(f)=\rho(\vee(x_j\iota+y_j))=\vee (\lambda x_j+y_j)=f(\lambda)$. \endproof
\begin{rem}{\rm One has in general, for a semiring $R$ of characteristic $1$, a natural isomorphism
 \begin{equation}\label{resiso}
{\rm Res}:\Hom_{\rma}(R\hat\otimes_\B\rma,\rma)\simeq \Hom_\B(R,\rma)
 \end{equation}
given by restriction of $\chi\in \Hom_{\rma}(R\hat\otimes_\B\rma,\rma)$ to the canonical image of $R$. 
Taking $R=\zmax$, this shows that the space $\Hom_{\rma}(\zmax\hat\otimes_\B\rma,\rma)$ of characters of $ \Z_{\rm max}\hat\otimes_\B\R_{\rm max}$   is the same as  $\Hom_\B(\zmax,\rma)$.}	
\end{rem}

\subsection{The small category $\mathscr C$}\label{sectcatC}

 The topos $\rnt$ is defined by assigning a small category $\mathscr C$ endowed with a Grothendieck topology $J$. We first describe $\mathscr C$. The objects of $\mathscr C$ are the (possibly empty) bounded open intervals 
$\Omega\subset [0,\infty)$ including those of the form $[0,a)$ for $a>0$. The morphisms between two objects of $\mathscr C$ are defined by 
$$
\Hom_{\mathscr C}(\Omega,\Omega')=\{n\in \nt\mid n\Omega\subset \Omega'\}
$$
 if $\Omega\neq \emptyset$. By definition  $\Hom_{\mathscr C}(\emptyset,\Omega'):=\{*\}$ \ie the one point set, for any object $\Omega'$ of $\mathscr C$. Thus the empty set is the initial object of $\mathscr C$.
The following lemma shows that pullbacks exist in the category $\mathscr C$
 \begin{lem}\label{lemcatC} 
 Let $\Omega_j\neq \emptyset$ ($j=1,2$) and consider two morphisms $\phi_j:\Omega_j\stackrel{n_j}{\to} \Omega$ given by integers $n_j\in \Hom_{\mathscr C}(\Omega_j,\Omega)$. Let $n={\rm lcm}(n_j)$ be their lowest common multiple,  $n=a_jn_j$ and let $\Omega':=\{\lambda \in [0,\infty)\mid a_j\lambda \in \Omega_j, \ j=1,2\}$. Then  $\Omega'$  is an object of $\mathscr C$,  and if it is non-empty one has $a_j\in \Hom_{\mathscr C}(\Omega',\Omega_j)$ and $(\Omega',a_j)$ is the pullback of the $\phi_j$. When $\Omega'=\emptyset$ the pullback of the $\phi_j$ is the initial object of $\mathscr C$. 
 \end{lem}
 \proof By construction $\Omega'=\cap_{j=1}^2 a_j^{-1}\Omega_j$ is an intersection of bounded open intervals and is thus an object of $\mathscr C$. Let $\psi_j:W\stackrel{k_j}{\to} \Omega_j$ be morphisms such that $\phi_1\circ \psi_1=\phi_2\circ \psi_2$. If $W\neq \emptyset$ it contains a $\lambda\neq 0$ and one has $n_1k_1=n_2k_2=mn$ for a unique $m\in \nt$. One has $k_j=ma_j$ and $k_jW\subset \Omega_j$ so that $mW\subset \Omega'$. Moreover the map $\psi_j:W\stackrel{k_j}{\to} \Omega_j$ is the composite of $W\stackrel{m}{\to} \Omega'$ with $a_j\in \Hom_{\mathscr C}(\Omega',\Omega_j)$.
 This shows that $(\Omega',a_j)$ is the pullback of the $\phi_j$. It also shows that if there exists $W\neq \emptyset$ and morphisms $\psi_j:W\stackrel{k_j}{\to} \Omega_j$
such that $\phi_1\circ \psi_1=\phi_2\circ \psi_2$, then $\Omega'\neq \emptyset$. Otherwise \ie if this implies $W=\emptyset$ then one easily sees that the empty set is indeed the pullback. \endproof
\subsection{The topology $J$ on $\mathscr C$}\label{sectGtop}
A Grothendieck topology $J$ (\cff \cite{MM} Definition III.2.1) on a small category $\mathscr C$ associates to every object 
$\Omega$ of the category a collection $J(\Omega)$ of sieves of $\Omega$, (\ie of families, stable under right composition, of morphisms with codomain $\Omega$),  such that:
 \vspace{.05in} 

$\blacktriangleright$~The maximal sieve $\{f\mid \Codom f=\Omega\}$ belongs to $J(\Omega)$
 \vspace{.05in} 

$\blacktriangleright$~$S\in J(\Omega)$, $h\in \Hom_{\mathscr C}(\Omega',\Omega)$ $\Rightarrow$ $h^*(S)\in J(\Omega')$, where
$
h^*(S):=\{f\mid h\circ f\in S\}
$ 

$\blacktriangleright$~For $S\in J(\Omega)$, and any sieve $R$ of $\Omega$ 
$$
h^*(R)\in J(\Dom\, h), \ \forall h\in S\ \Rightarrow R\in J(\Omega).
$$
When the small category $\mathscr C$ admits pullbacks one can associate a Grothendieck topology $J$ to a basis $K$, \ie a function  which assigns to any object $\Omega$ a collection $K(\Omega)$ of families of morphisms with codomain $\Omega$ by the condition 
$$
S\in J(\Omega) \iff \exists R\in K(\Omega), \ \ R\subset S.
$$
The above three conditions on $J$ are derived from the following three conditions on $K$: (\cff \cite{MM} Definition III.2.2):
\begin{enumerate}
\item For any isomorphism $f$ with range $\Omega$ the singleton $\{f\}$ is a covering. 
\item The pullback of a covering of $\Omega$ by any morphism $\Omega'\to\Omega$ is a covering of $\Omega'$.
\item Given a covering $(\Omega_j)_{j\in I}$ of $\Omega$ and for each $j\in I$ a covering $\Omega_{ij}$ of $\Omega_j$, the composite $\Omega_{ij}$ is a covering of $\Omega$.
\end{enumerate}
\begin{prop}\label{proptop} 
$(i)$~For each object $\Omega$ of $\mathscr C$, let $K(\Omega)$ be the collection of all ordinary covers $\{\Omega_i\subset \Omega, i\in I\mid \cup \Omega_i=\Omega\}$ of $\Omega$. Then $K$ defines a basis for a Grothendieck topology $J$ on $\mathscr C$.\newline
$(ii)$~The Grothendieck topology $J$ is subcanonical.\newline
$(iii)$~The category $\sh(\mathscr C,J)$ of sheaves of sets on $(\mathscr C,J)$ is canonically isomorphic to the category of $\nt$-equivariant sheaves of sets on $[0,\infty)$.	
\end{prop}
\proof $(i)$~The only isomorphisms in $\mathscr C$ are the identity maps, thus one verifies 1. To check the condition 2., we let $\Omega' \stackrel{n}{\to} \Omega$ be a morphism in $\mathscr C$ and 
$\{\Omega_i\subset \Omega, i\in I\mid \cup \Omega_i=\Omega\}$ a covering of $\Omega$. Then it follows from Lemma \ref{lemcatC} that the pullback of the cover is given by 
$$
\left(\pi_2:\Omega_i\times_\Omega \Omega'\to \Omega' \right)\simeq \left(n^{-1}\Omega_i\cap \Omega'\to \Omega' \right).
$$
This defines a covering of $\Omega'$. Finally the condition 3. is a standard fact on ordinary covers of a topological space (here chosen to be $[0,\infty)$).\newline
$(ii)$~We prove that any representable presheaf is a sheaf on $(\mathscr C,J)$. For a fixed object $\Omega$ of $\mathscr C$ and an arbitrary open subset $U$ of $[0,\infty)$ one sets
$$
\Gamma(U):=\Hom_{\mathscr C}(U,\Omega)=\{n\in \nt\mid nU\subset \Omega\}.
$$
 This determines the subsheaf of the constant sheaf $\nt$ which is given by the local condition around $\lambda$: $\{n\in \nt\mid n\lambda \in \Omega\}$. \newline
 $(iii)$~An $\nt$-equivariant sheaf (of sets) on $[0,\infty)$ gives by restriction to $\mathscr C$ an object of $\sh(\mathscr C,J)$. Conversely let $\cF$ be an object of $\sh(\mathscr C,J)$, and $U$ an arbitrary open subset of $[0,\infty)$. Let $\{\Omega_i\subset U, i\in I\mid \cup \Omega_i=U\}$ be a covering of $U$ by bounded open intervals. Take
 the limit $\varprojlim \cF(\Omega_i\cap\Omega_j)$ in $\Se$ of the diagram $\cF(\Omega_i\cap\Omega_j)$ indexed by pairs $(i,j)\in I^2$, $\Omega_i\cap \Omega_j\neq \emptyset$ and arrows $(i,j)\to (i,i)$, $(i,j)\to (j,j)$. This is the equalizer of the two maps
 $$
 \prod_{i\in I} \cF(\Omega_i)\stackrel{p_\ell}{\rightrightarrows}\prod_{i,j}\cF(\Omega_i\cap\Omega_j).
 $$
 Since $\cF$ is an object of $\sh(\mathscr C,J)$ the above equalizer does not depend upon the choice of the covering of $U$ by bounded open intervals, and  defines a sheaf of sets on $[0,\infty)$ endowed with an action of $\nt$ compatible with the action of $\nt$ on $[0,\infty)$.\endproof
 
 \begin{rem}\label{empty3}{\rm The setting $\Hom_{\mathscr C}(\emptyset,\Omega)= \{*\}$ is imposed if one wants the Grothendieck topology $J$  to be subcanonical. Indeed, any sheaf evaluated on the empty set gives the one point set $\{*\}$. 
 }	
 \end{rem}

\section{The points of the topos $\rnt$} \label{sectptsss}

In this section we investigate the structure of the points of the  topos $\rnt$ and prove Theorem \ref{scaltop}
which provides a canonical bijection between the space of these points (up to isomorphism) and  the sector $\Q^\times\backslash\A_\Q/\hat\Z^*$  of the adele class space of $\Q$.

\subsection{Flatness and continuity}
It follows from \cite{MM} Corollary VII.5.4, that the points of the topos $\sh(\mathscr C,J)$ correspond by means of an equivalence of categories to continuous, flat functors $F:\mathscr C\longrightarrow \Se$. Moreover again from \cite{MM} Theorem VII.6.3, one knows that a functor $F:\mathscr C\longrightarrow \Se$ is flat iff it is filtering, \ie $F$ fulfills the  three conditions reported in the following definition
\begin{defn} \label{defnfiltering} A functor $F:\mathscr C\longrightarrow \Se$ is filtering iff 
\begin{enumerate}
	\item $F(C)\neq \emptyset$ for some object $C$ of $\mathscr C$. 
	\item Given $a_j\in F(C_j)$ $j=1,2$, there exist an object $C$ of $\mathscr C$, an element $a\in F(C)$  and morphisms $u_j:C\to C_j$ such that $F(u_j)a=a_j$.
	\item Given two morphisms $u,v:C\to D$ in $\mathscr C$ and $a\in F(C)$ such that $F(u)a=F(v)a$, there exists an object $B$ of $\mathscr C$, an element $b\in F(B)$ and a morphism $w:B\to C$ of $\mathscr C$ such that $u\circ w=v\circ w$ and $F(w)b=a$.
\end{enumerate}
\end{defn}
By \cite{MM}, Lemma VII.5.3, a flat functor $F:\mathscr C\longrightarrow \Se$ is continuous iff it sends covering sieves to epimorphic families. 
\begin{lem}\label{lemcont1} Let $F:\mathscr C\longrightarrow \Se$ be a flat functor and let $J$ be the Grothendieck topology on $\mathscr C$ generated by a basis $K$. Then $F$ is continuous iff for any object $U$ of $\mathscr C$ and a covering $(U_j\to U)_{j\in I}\in K(U)$ the family of maps $F(U_j\to U):F(U_j)\to F(U)$ is jointly surjective.	
\end{lem}
\proof Since any covering sieve $S\in J(U)$ contains a covering from the basis, the condition of the Lemma implies that $F$ sends covering sieves to epimorphic families. Conversely let $R\in K(U)$, then the associated covering sieve $(R)\in J(U)$
$$
(R):=\{f\circ g\mid f\in R, \ \Dom \, f=\Codom\,  g\}
$$
gives an epimorphic family $F(f\circ g)$ iff the family $F(f)$ is jointly surjective.\endproof

\subsection{The point $\ffp_H$ associated to a rank one subgroup of $\R$}

The next Proposition shows that any (non-trivial) rank one subgroup $H\subset \R$ defines a point of the topos $\rnt$.

\begin{prop} \label{proppoint} 
$(i)$~Let $H$ be a (non-trivial) rank one subgroup of $\R$, then the equality $F_H(V):=V\cap H \cap (0,\infty)$ defines a flat continuous functor $F_H: \mathscr C\longrightarrow \Se$.

$(ii)$~The map $H\mapsto \ffp_H$ which associates to a rank one subgroup of $\R$ the point of the topos $\rnt$ represented by the flat continuous functor $F_H$ provides an injection of the space of (non-trivial) rank one subgroups of $\R$ into the space of points of the topos $\rnt$ up to isomorphism.
\end{prop}
\proof $(i)$~We set $H_+:=H \cap (0,\infty)$. Let $V \stackrel{n}{\to} W$ be a morphism in $\mathscr C$, then since $n H_+\subset H_+$ the equality $F_H(V)=V\cap H_+$ defines a subfunctor of the covariant functor $\mathscr C\longrightarrow \Se$, $V\mapsto V$. Given a covering $(U_j\subset U)_{j\in I}\in K(U)$ the family of maps $F_H(U_j\to U):F_H(U_j)\to F_H(U)$ is jointly surjective since $U\cap H_+=\cup(U_j\cap H_+)$, thus Lemma \ref{lemcont1} shows that $F_H$ is continuous. We show that $F_H$ is filtering. Since $H$ is non-trivial, $H_+\neq \emptyset$ and this gives condition $1$. Next, given $h_j\in V_j\cap H_+$ we let $h\in H_+$, $n_j\in \nt$ be such that $h_j=n_j h$. Let $V$ be  an open interval containing $h$ and such that $n_j V\subset V_j$. This defines an object $V$ of $\mathscr C$. One has $h\in V\cap H_+=F_H(V)$ and the morphisms $u_j:V \stackrel{n_j}{\to} V_j$ fulfill $F_H(u_j)h=h_j$. This gives condition $2$. Finally, if  $u:C \stackrel{n}{\to} D$, $v:C \stackrel{m}{\to} D$ are two morphisms in $\mathscr C$ and $a\in F_H(C)$ is such that $F_H(u)a=F_H(v)a$, then since $a>0$ one gets that $n=m$ and thus the condition $3$ holds. \newline
$(ii)$~Let $F_H(V):=V\cap H_+$ be the continuous flat functor associated to $H\subset \R$. Given a point $\lambda\in (0,\infty)$, we let $V_j$ be a basis of neighborhoods of $\lambda$ of bounded open intervals. Then one has
$$
\varprojlim_{U\ni x} F_H(U)\neq \emptyset \iff\cap F_H(V_j)\neq \emptyset \iff \lambda \in H.
$$
This shows that one can recover the subgroup $H\subset \R$ from the continuous flat functor $F_H$. Moreover it shows that a morphism of functors from $F_H$ to $F_{H'}$ exists iff $H\subset H'$ and hence that  the isomorphism class of the point $\ffp_H$ uniquely determines $H\subset \R$.\endproof

\subsection{Classification of points of the topos $\rnt$}
The main result of this subsection (Theorem~\ref{scaltop}) states the existence  of a canonical isomorphism between the points $\mathscr A(\rmax)$ of the arithmetic site over $\rmax$ and the points of the topos $\rnt$. To this end we shall need to state first several technical lemmas.
\begin{lem}\label{lemflatcont1}Let $F:\mathscr C\longrightarrow \Se$ be a continuous flat functor. Then the following facts hold\newline
$(i)$~Let $U_j\subset V$ ($j=1,2$) be objects of $\mathscr C$ with $U_1\cap U_2=\emptyset$. Then the images $F(U_j\to V)F(U_j)\subset F(V)$ are disjoint.\newline
$(ii)$~Let $V$ be an object of $\mathscr C$ and $x\in F(V)$. Then there exist a unique $\lambda\in V$ such that for any object $U$ of $\mathscr C$, $U\subset V$,  containing $\lambda$ one has $x\in F(U\to V)F(U)\subset F(V)$.
\end{lem}
\proof $(i)$~Since $F$ is flat it commutes with fibered products. The fibered product of the two maps $U_j\to V$ is the empty set. The object $\emptyset$ of $\mathscr C$ admits the empty cover. Thus by continuity one has $F(\emptyset)=\emptyset$. One can also give the following argument. Let $a_j\in U_j$ with $ F(U_1\to V)(a_1)=F(U_2\to V)(a_2)=z$. Then the flatness of $F$ gives an object $W$ of $\mathscr C$ an element $c\in F(W)$ and morphisms $W \stackrel{n_j}{\to} U_j$ such that $F(W \stackrel{n_j}{\to} U_j)(c)=a_j$. By composition with $U_j\to V$ one gets $F(W \stackrel{n_j}{\to} V)(c)=z$. The third filtering condition on $F$ implies, since $\nt$ is simplifiable, that $n_1=n_2$ which contradicts $n_jW\subset U_j$ and $U_1\cap U_2=\emptyset$. \newline
$(ii)$~We show the existence of $\lambda$, its uniqueness then follows from $(i)$ since distinct points have disjoint neighborhoods given by objects of $\mathscr C$. Let  $V=\cup W_j$ where $W_j$ is an increasing family of bounded open intervals such that $\overline{ W_j}\subset W_{j+1}$. Then the continuity of $F$ gives an interval $W=W_j$, with $\overline{W}\subset V$ such that $x\in F(W\to V)F(W)$, \ie $x=F(W\to V)z$, $z\in F(W)$. Using a cover $\cU$ of $W$ by bounded open intervals   one obtains, by continuity of $F$,  an interval $I_1\subset W$ of diameter $<1/2$ and $z_1\in F(I_1)$ such that $z= F(I_1\to W)z_1$. By induction one gets a decreasing sequence of intervals $I_k\subset I_{k-1}\subset W$ of diameter $<1/2^k$ and $z_k\in F(I_k)$ such that $z_{k-1}= F(I_k\to I_{k-1})z_k$. Let $\lambda$ be  the unique limit point of the sequence $I_k$, \ie the limit of any sequence $\lambda_k\in I_k$. One has $\lambda\in\overline{W}\subset V$. Let $U\subset V$ be an object of $\mathscr C$ containing $\lambda$. Then $U$ is an open neighborhood of $\lambda$ and  there exists  $k$ such that $I_k\subset U$. One then gets
$$
x=F(W\to V)z=F(I_k\to V)z_k=F(U\to V)F(I_k\to U)z_k\in F(U\to V)F(U)
$$ 
\endproof
In what  follows we shall  denote by $\lambda_V:F(V)\to V$ the map defined in Lemma \ref{lemflatcont1}.
\begin{lem}\label{lemflatcont2}Let $F:\mathscr C\longrightarrow \Se$ be a continuous flat functor. \newline
$(i)$~The maps $\lambda_V:F(V)\to V$ define a natural transformation of $F$ with the functor $\mathscr C\longrightarrow \Se$, $V\mapsto V$.\newline
$(ii)$~The maps $\lambda_V:F(V)\to V$ are injective when $0\notin V$.
\end{lem}
\proof $(i)$~Let $U\stackrel{n}{\to} V$ be a morphism in $\mathscr C$. Let $x\in F(U)$, $y=F(U\stackrel{n}{\to} V)x\in F(V)$. For any object $W\subset U$ of $\mathscr C$ containing $\lambda_U(x)$ one has $x\in F(W\to U)F(W)$ and thus 
$$
y\in F(U\stackrel{n}{\to} V)F(W\to U)F(W)= F(W\stackrel{n}{\to} V)F(W)
=F(nW\to V)F(W\stackrel{n}{\to} nW)F(W).
$$
This shows that $y\in F(nW\to V)F(nW)$. Thus if $\lambda_V(y)\neq n \lambda_U(x)$ one obtains a contradiction by Lemma \ref{lemflatcont1} $(i)$.
We have thus shown that the following diagram commutes
\begin{equation*}\label{1o}
\xymatrix{
F(U) \ar[d]^{\lambda_U}  \ar[rr]^{F(U\stackrel{n}{\to} V)}&& \ar[d]^{\lambda_V}  F(V) \\
U \ar[rr]^{n}&& V
}
\end{equation*}
$(ii)$~Let $x_j\in F(V)$. By the flatness of $F$ there exists an object $W$ of $\mathscr C$, an element $c\in F(W)$ and morphisms $W \stackrel{n_j}{\to} V$ such that $F(W \stackrel{n_j}{\to} V)(c)=x_j$. Since $0\notin V$ by hypothesis, one has $\lambda_V(x_1)\neq 0$. Assume now that $\lambda_V(x_1)=\lambda_V(x_2)$. By $(i)$ one has $\lambda_V(x_j)=n_j\lambda_W(c)$ and this implies $\lambda_W(c)\neq 0$ (since $\lambda_V(x_1)\neq 0$) and $n_1=n_2$. One then gets $x_1=x_2$ since one has $x_j=F(W \stackrel{n_j}{\to} V)(c)$.\endproof 
\begin{lem}\label{lemflatcont3}Let $F:\mathscr C\longrightarrow \Se$ be a continuous flat functor. \newline
$(i)$~Let $\lambda>0$ be a positive real number then, for objects $V$ of $\mathscr C$ one has 
$$
\exists V\mid \lambda\in \lambda_V(F(V))\iff \forall V\ni \lambda, \ \lambda\in \lambda_V(F(V))
$$
$(ii)$~The subset of $(0,\infty)$ defined by the above condition is of the form $H\cap (0,\infty)$, where $H\subset \R$ is a rank one subgroup.
\end{lem}
\proof $(i)$~It is enough to show the implication $\Rightarrow$. Let $V$ and $x\in F(V)$ such that $\lambda=\lambda_V(x)$. Then for any object $U$ of $\mathscr C$, $U\subset V$,  containing $\lambda$ one has $x\in F(U\to V)F(U)\subset F(V)$. Let then $W$ be an object of $\mathscr C$ containing $\lambda$. Let $U=W\cap V$ and $z\in F(U)$ such that $x=F(U\to V)z$. One has by Lemma \ref{lemflatcont2} $(i)$, $\lambda_W(F(U\to W)z)=\lambda_U(z)$ and $\lambda_U(z)=\lambda_V(F(U\to V)z)=\lambda_V(x)=\lambda$. Thus one gets $\lambda\in \lambda_W(F(W))$ as required. 

$(ii)$~Let $E=\{\lambda>0\mid \lambda\in \lambda_V(F(V))\qqq V\ni \lambda\}$. Let $\lambda_j\in E$, $j=1,2$, and $V$ containing the $\lambda_j$, and 
$x_j\in F(V)$ such that $\lambda_j=\lambda_V(x_j)$. By the flatness of $F$ there exists an object $W$ of $\mathscr C$, an element $c\in F(W)$ and morphisms $W \stackrel{n_j}{\to} V$ such that $F(W \stackrel{n_j}{\to} V)(c)=x_j$.  By Lemma \ref{lemflatcont2} $(i)$, one has $\lambda_V(x_j)=n_j\lambda_W(c)$. By $(i)$ one has $\lambda_W(c)\in E$. This shows   that given any two elements  $\lambda_j\in E$ there exists $\lambda \in E$ and integers $n_j\in \nt$ such that $n_j\lambda =\lambda_j$, $j=1,2$. Moreover Lemma \ref{lemflatcont2} $(i)$ shows that $\lambda\in E\Rightarrow n\lambda \in E$, $\forall n\in \nt$. Thus $E$ is an increasing union of subsets of the form $\lambda_k \nt$. Let 
$H=\cup \lambda_k\Z$ be the corresponding increasing union of subgroups of $\R$. Then $H$ is a rank one subgroup of $\R$ and $E=(0,\infty)\cap H$ by construction. \endproof

\begin{lem}\label{lemflatcont4}Let $F:\mathscr C\longrightarrow \Se$ be a continuous flat functor. \newline
$(i)$~One has 
$$
\exists V\mid 0\in \lambda_V(F(V))\iff \forall V, \ \lambda_V(F(V))=V\cap \{0\}
$$
$(ii)$~If the above equivalent conditions do not hold then there exists a rank one subgroup $H\subset \R$ and an isomorphism of functors $F\simeq F_H$.
\end{lem}
\proof $(i)$~Let $V$ be such that $0\in \lambda_V(F(V))$. One then has $0\in V$. Moreover the proof of Lemma \ref{lemflatcont3} $(i)$ applies and shows that for any object $W$ of $\mathscr C$ containing $0$ one has $0\in \lambda_W(F(W))$. Let $W$ be an object of $\mathscr C$  and assume that some $\lambda>0$ belongs to 
$\lambda_W(F(W))$. Let $U$ be an object of $\mathscr C$ containing both $V$ and $W$. Then by Lemma \ref{lemflatcont2} $(i)$, one has $\{0,\lambda\}\subset \lambda_U(F(U))$. Let $x_j\in F(U)$ with $\lambda_U(x_1)=0$, $\lambda_U(x_2)=\lambda$. Then the flatness of $F$ gives an object $U'$ of $\mathscr C$, an element $c\in F(U')$ and morphisms $U' \stackrel{n_j}{\to} U$ such that $F(U'\stackrel{n_j}{\to} U)(c)=x_j$.  By Lemma \ref{lemflatcont2} $(i)$, one has $\lambda_U(x_j)=n_j\lambda_{U'}(c)$. But $\lambda_U(x_1)=0$  implies that $\lambda_{U'}(c)=0$ and this contradicts $\lambda_U(x_2)=\lambda\neq 0$. Thus it follows that for any object $W$ of $\mathscr C$, $\lambda_W(F(W))$ contains at most $0$ and it does if and only if one has $0\in W$. \newline
 $(ii)$~If the condition of $(i)$ does not hold, it follows that for any object $V$ of $\mathscr C$ one has $0\notin \lambda_V(F(V))$. It follows that the canonical map $\rho_V=F(V\cap (0,\infty)\to V): F(V\cap (0,\infty))\to F(V)$ is surjective. But by Lemma \ref{lemflatcont2}, $(ii)$, and the commutation with the localization $\lambda_V$ this map is injective. By Lemma \ref{lemflatcont3} $(ii)$ there exists a rank one subgroup $H\subset \R$ such that for any $U$ not containing $0$ one has $\lambda_U(F(U))=U\cap H$. By Lemma \ref{lemflatcont2}, $(ii)$, the composite 
 $$
 \lambda_{V\cap (0,\infty)}\circ \rho_V^{-1}: F(V)\to V\cap (0,\infty)\cap H
 $$
 gives an isomorphism of functors $F\simeq F_H$. \endproof 
 
 \begin{lem}\label{lemflatcont5}Let $F:\mathscr C\longrightarrow \Se$ be a continuous flat functor such that $\lambda_V(F(V))=V\cap \{0\}$, $\forall V$. \newline
$(i)$~One has $F(V)=\emptyset$ if $0\notin V$ and if $0\in V$ then the canonical map 
$p_V$ from $X:=\varprojlim_{W\ni 0} F(W)$ to $F(V)$ is bijective.\newline
$(ii)$~There exists a unique action $n\mapsto X(n)$ of $\nt$ on $X$ such that for any morphism $U \stackrel{n}{\to} V$ of objects containing $0$ the following diagram commutes :
\begin{equation}\label{1o1}
\xymatrix{
 X\ar[d]^{p_U}  \ar[rr]^{X(n)}&& \ar[d]^{p_V}  X \\
F(U) \ar[rr]^{F(U\stackrel{n}{\to} V)}&& F(V)
}
\end{equation}
$(iii)$~With the above notations, the action of $\nt$ on $X$ defines a point of the topos $\wnt$.\newline
$(iv)$~Let $H$ be an abstract rank one ordered group. The following defines a flat continuous functor $F'_H:\mathscr C\longrightarrow \Se$
$$
F'_H(V)=\emptyset \ \text{if} \ 0\notin V, \  \ F'_H(V)=H_+ \ \text{if} \ 0\in V.
$$
\end{lem}
\proof $(i)$~If $0\notin V$ one has $\lambda_V(F(V))=\emptyset$ and hence $F(V)=\emptyset$. Let $U\subset V$ with $0\in U$. The map $F(U\to V):F(U)\to F(V)$ is surjective since $\lambda_V(F(V))=\{0\}$. Let us show that it is injective. Let $x_j\in F(U)$ be such that $F(U\to V)(x_j)=z$. By the flatness of $F$ there exists an object $W$ of $\mathscr C$, an element $c\in F(W)$ and morphisms $W \stackrel{n_j}{\to} U$ such that $F(W \stackrel{n_j}{\to} U)(c)=x_j$. One then has $F(W \stackrel{n_j}{\to} V)(c)=z$ and the flatness of $F$ shows, since $\nt$ is simplifiable, that $n_1=n_2$ so that $x_1=x_2$. Thus the map $F(U\to V):F(U)\to F(V)$ is bijective and the projective limit
$X:=\varprojlim_{W\ni 0} F(W)$ is such that all maps $p_V: X\to F(V)$ are bijective.  \newline
$(ii)$~The inclusions $U\subset V$, $nU\subset nV$ form a commutative square  with the maps $U \stackrel{n}{\to} nU$, $V \stackrel{n}{\to} nV$. It follows that the map from $X$ to $X$ such that 
$$
F(U \stackrel{n}{\to} nU)p_U(x)=p_{nU}(X(n)x)
$$
is independent of the choice of $U$ containing $0$ and turns \eqref{1o1} into commutative diagrams. \newline
$(iii)$~The set $X$ is non-empty since otherwise one would have $F(V)=\emptyset$, $\forall V$. Let $x_j\in X$ and $u_j\in F(U)$ ($0\in U$) such that $p_U(x_j)=u_j$.  By the flatness of $F$ there exists an object $W$ of $\mathscr C$, an element $c\in F(W)$ and morphisms $W \stackrel{n_j}{\to} U$ such that $F(W \stackrel{n_j}{\to} U)(c)=u_j$. One has $0\in W$. Let $x\in X$ such that $p_W(x)=c$, using \eqref{1o1} one gets $X(n_j)x=x_j$. Thus the action of $\nt$ on $X$ verifies the second filtering condition. We now check the third filtering condition.  Let $x\in X$ and $n_1\neq n_2$ be such that $X(n_1)x=X(n_2)x$. Let $0\in U$, $u=p_U(x)$. Let $V\supset n_jU$. One then gets  an equality of the form 
$$
F(U \stackrel{n_1}{\to} V)(u)=F(U \stackrel{n_2}{\to} V)(u)
$$
and the flatness of $F$ shows that $n_1=n_2$. This shows that the action of $\nt$ on $X$ verifies the three filtering conditions.\newline
$(iv)$~One checks that $F'_H$ is continuous, its flatness follows from the rank one property of $H$. \endproof 

We are now ready to define a map $\Theta$ from points of $\aarith(\rmax)$ to points of the  topos $\rnt$. As shown in \cite{CC1} Theorem 3.8,  the points of the arithmetic site  $\aarith$ over $\rmax$ form the union of two sets: the set $\aarith(\B)\subset \aarith(\rmax)$ of isomorphism classes of points of $\wnt$,  and the set of the non-trivial rank one subgroups $H\subset \R$. To a point of $\aarith(\B)\subset \aarith(\rmax)$ associated with an abstract rank one ordered group $H$ we assign the point $\ffq_H$ of $\rnt$ associated to the flat continuous functor $F'_H$ of 
 Lemma \ref{lemflatcont5}, $(iv)$. 
Next to the point of $ \aarith(\rmax)\setminus \aarith(\B)$ corresponding to the rank one subgroup $H\subset \R$, we associate the point $\ffp_H$ of Proposition \ref{proppoint}.
\begin{thm}\label{scaltop}The  map  $\Theta$ defines a canonical isomorphism of $\aarith(\rmax)$ with the points of the  topos $\rnt$.
\end{thm}
\proof The map $\Theta$ is well defined. Lemma \ref{lemflatcont5} shows that is bijective from points of $\aarith(\B)\subset \aarith(\rmax)$ to 
points of $\rnt$ such that the associated flat continuous functor $F:\mathscr C\longrightarrow \Se$ fulfills the hypothesis of Lemma \ref{lemflatcont5}. Proposition \ref{proppoint}  shows that $\Theta$ is injective from points of $\aarith(\rmax)\setminus \aarith(\B)$ to 
points of $\rnt$. Finally, Lemma \ref{lemflatcont4} shows that $\Theta$  is surjective.\endproof

\subsection{Representation of points as filtering colimits of representables}\label{sectrep}

Let $\mathscr C$ be a small category. Any  covariant functor $F:\mathscr C\longrightarrow \Se$ which is 
representable, \ie of the form $y_I(C)=\Hom_\mathscr C(I,C)$ for some object $I$ of $\mathscr C$ is flat, and any flat functor is obtained as a filtering colimit of such representable functors. In this subsection we describe such representations for the flat continuous functors associated to points of the topos $\rnt$. This result will be used in the description of the stalks of the sheaves given in \S~\ref{sectsheaf}. 
\begin{lem} \label{limlim1} 
Let $H\subset \R$ be a rank one subgroup of $\R$. Let $h_i\in H_+$ be a sequence of elements, $n_j\in\nt$, such that $n_jh_{j+1}= h_j$ and that $H=\cup h_j\Z$. Let $I_j$ be bounded open intervals with $h_j\in I_j$ for $j\geq 1$, such that 
\begin{equation}\label{gofast}
n_j\overline{ I_{j+1}}\subset I_j, \ \forall j\geq 1, \  \ 
\,(\prod_1^{k-1} n_i)\, {\rm Diameter}(I_k)\to 0, \  \ \text{when}\, \ k\to \infty.
\end{equation}
Then the limit $\varinjlim y_{I_j}$ of the representable functors $y_{I_j}(V):=\Hom_\mathscr C(I_j,V)$ defines the point 
$\ffp_H$ of the topos $\rnt$.
\end{lem}
\proof We show that the functor $\varinjlim y_{I_j}$ is simply given by 
$$
\varinjlim y_{I_j}(V)=\varinjlim( h_j^{-1} V\cap \nt)\sim \varinjlim(  V\cap h_j \nt)=V\cap H_+.
$$
It is enough to prove that the natural inclusion, due to $h_j\in I_j$ 
$$
\varinjlim y_{I_j}(V)\subset \varinjlim( h_j^{-1} V\cap \nt)
$$
(in both cases one uses multiplication by $n_j$ to organize the inductive system) is a bijection. 
Indeed, let $n\in h_k^{-1} V\cap \nt$: we show that for $j$ large enough one has $n\prod_k^{j-1} n_i\in y_{I_j}(V)$. Let $\epsilon>0$ be a positive real number such that the neighborhood $W$ of $h_k n$ of radius $\epsilon$ is contained in $V$. Then for $j$ large enough, using the hypothesis \eqref{gofast}, one gets that $n\prod_k^{j-1} n_i\, I_j\subset W$ (since it contains $h_kn$ and is of small enough diameter) and hence $n\prod_k^{j-1} n_i\in y_{I_j}(V)$. \endproof

\begin{prop}\label{proplim} Let $H\subset \R$ be a rank one subgroup of $\R$. Let $h_i\in H_+$, $n_j\in\nt$ and $I_j$ be bounded open intervals as in Lemma \ref{limlim1}. Then, the pullback part of the point $\ffp_H$ of the topos $\rnt$ is given by the functor which associates to an object $\cF$ of $\sh(\mathscr C,J)$ the following colimit
\begin{equation}\label{gofast1}
\varinjlim_k \cF(I_k), \quad  \cF(I_{k+1} \stackrel{n_k}{\to} I_k):\cF(I_k)\to \cF(I_{k+1}).
\end{equation}	
\end{prop}
\proof The colimit \eqref{gofast1} defines the pullback part $f^*$ of a point of the topos of contravariant functors $\mathscr C\longrightarrow \Se$, which is defined as a filtering colimit of the points associated to the objects $I_k$ of $\mathscr C$. To show that the corresponding geometric morphism from the topos of sets to $\hat{\mathscr C}$ factors through $\sh(\mathscr C,J)$ it is enough to show (\cf~\cite{MM} Lemma VII.5.3) that the composite $f^*\circ y$ with the Yoneda embedding sends each covering sieve to an epimorphic family of functions. Lemma \ref{limlim1} shows that for any object $V$ of $\mathscr C$ and the associated object of $\hat{\mathscr C}$: $\cF=y(V):=\Hom_\mathscr C(\bullet, V)$ one has 
$$
\varinjlim_k \cF(I_k)= \varinjlim_k \Hom_\mathscr C(I_k,V)=V\cap H_+.
$$
 Thus  $f$ fulfills the condition $(iii)$ of \cite{MM} Lemma VII.5.3, and  $f^*$ is the pullback part of the point $\ffp_H$ of the topos $\rnt$. \endproof
 Next, we consider the limit case of Lemma \ref{limlim1} when all the $h_j$   are $0$.
 \begin{lem} \label{limlim2} 
Let  $n_j\in\nt$, and $I_j$ be bounded open intervals containing $0\in[0,\infty)$,  for $j\geq 1$, such that \eqref{gofast} holds.
Then the limit $\varinjlim y_{I_j}$ of the representable functors $y_{I_j}(V):=\Hom_\mathscr C(I_j,V)$ defines the point $\ffq_H$ 
 of the topos $\rnt$ associated to the abstract ordered rank one group $H:=\cup (\prod_1^k n_j)^{-1}\Z$.
\end{lem} 
\proof Let $F=\varinjlim y_{I_j}$. One has  $F(V)=\emptyset$ when $0\notin V$. Let $m_j:=\prod_1^{j-1} n_\ell$ and $H:=\cup\, m_j^{-1}\Z$.  We denote by $\iota(j,k):=k/m_j\in H$ the element of $H$ given by  the image of $k\in \Z$ by the canonical injection $\Z\to H$ associated to the $j$-th copy of $\Z$ in the colimit. One has by construction 
\begin{equation}\label{iotak}
\iota(j,k)=\iota(j+1, n_j\, k)\qqq j, \ k\in \Z	.
\end{equation}
 Let $V$ be a bounded open interval with $0\in V$. One then obtains an injection 
 $$
 \alpha:F(V)\to H_+, \  \  \alpha(I_j\stackrel{k}{\to} V):=\iota(j,k)
 $$
by using the compatibility with the inductive limits, \ie 
$$
\alpha\left( (I_j\stackrel{k}{\to} V)\circ (I_{j+1}\stackrel{n_j}{\to} I_j) \right)
=\iota(j+1, n_j\, k)=\iota(j,k).
$$
Let us show that $\alpha$ is surjective. 
 Given $k\in \nt$ and $j>0$ we prove that for $\ell$ large enough one has $k\prod_j^{\ell-1} n_i\in y_{I_\ell}(V)$. Let $\epsilon>0$ such that the neighborhood $W$ of $0$ of radius $\epsilon$ is contained in $V$. Then for $\ell$ large enough, using the hypothesis \eqref{gofast}, one gets that $k\prod_j^{\ell-1} n_i\, I_\ell\subset W$ (since it contains $0$ and is of small enough diameter) and hence $k\prod_j^{\ell-1} n_i\in y_{I_\ell}(V)$. \endproof
 \begin{prop}\label{proplimbis} Let $H$ be an abstract rank one ordered group. Let  $n_j\in\nt$, $I_j$ be bounded open intervals as in Lemma \ref{limlim2}  such that $H:=\cup (\prod_1^k n_j)^{-1}\Z$. Then, the pullback part of the point $\ffq_H$ of the topos $\rnt$ is given by the functor which associates to an object $\cF$ of $\sh(\mathscr C,J)$ the colimit \eqref{gofast1}.
\end{prop}
The proof is the same as the one of Proposition \ref{proplim}. 

\section{The structure sheaf $\cO$ of the scaling site} \label{sectsheaf}

We define the structure sheaf $\cO$ of the scaling site $\scal2$ as the $\nt$-equivariant sheaf on  $[0,\infty)$ of semirings $\cO(U)$ of continuous convex   functions  $f(\lambda)\in \rma$, such that  for any $\lambda \in U$ there exists an open interval $V$ containing $\lambda$, with  $f$  affine and with slope $f'\in \Z$ in the complement  $V\setminus \{\lambda\}$.
This condition is local and hence defines a sheaf. One endows this sheaf  with the following action $\cO(V\stackrel{n}{\to}W): \cO(W)\to \cO(V)$ of $\nt$
\begin{equation}\label{equivO}
	\cO(V\stackrel{n}{\to}W)(f)(\lambda):= f(n\lambda)\qqq \lambda \in V.
\end{equation}
This action is compatible with the semiring structure and with the integrality property  of the slopes and thus defines a semi-ring in the topos $\rnt$ 
\begin{defn}\label{defnscalsite} The scaling site $\scal2$ is the semi-ringed topos given by the topos $\rnt$ endowed with the structure sheaf $\cO$.	
\end{defn}

\subsection{The stalks of the structure sheaf $\cO$}\label{sectsheafstalks}
The next result determines the structure of the stalks of $\cO$.

\begin{thm} \label{structure2} $(i)$~At the point $\ffp_H$ of the topos $\rnt$ associated to the rank one subgroup $H\subset \R$ the stalk of the structure sheaf $\cO$  is the semiring $\cR_H$ of germs, at $\lambda=1$, of $\rma$-valued, piecewise affine   convex functions $f(\lambda)$ with slopes in $H$.\newline
$(ii)$~The stalk of the structure sheaf $\cO$ at the point $\ffq_H$ of $\rnt$ associated to the abstract rank one ordered group $H$ is the semiring $\cZ_H$ associated by the max-plus construction to the totally ordered group $\R\times H$ and endowed with the lexicographic order.
\end{thm}
\proof $(i)$~To evaluate the stalk of the structure sheaf $\cO$ at the point $\ffp_H$ we use the description  given by \eqref{gofast1}. We let $h_j$, $n_j$, $I_j$ as in Proposition \ref{proplim} and evaluate the colimit
\begin{equation*}\label{gofast2}
\cO_{\ffp_H}=\varinjlim_k \cO(I_k), \  \  \cO(I_{k+1} \stackrel{n_k}{\to} I_k):\cO(I_k)\to \cO(I_{k+1}).
\end{equation*}
We define a map $\rho:\cO_{\ffp_H}\to \cR_H$ by associating to $(j,f)$, $f\in \cO(I_j)$, the germ at $\lambda=1$ of the function $\lambda\mapsto f(\lambda h_j)$. This function is defined in the neighborhood of $\lambda=1$ given by $\{\lambda\mid h_j\lambda \in I_j\}$, and it is a piecewise affine, continuous, convex function with slopes in $h_j\Z\subset H$. Thus its germ  at $\lambda=1$  is an element  $\rho(j,f)\in \cR_H$. Next we prove that this construction is compatible with the colimit, \ie that 
$$
\rho(j,f)=\rho(j+1,\cO(I_{j+1} \stackrel{n_j}{\to} I_j)(f))
\qqq f\in \cO(I_j)$$ 
where $\cO(I_{j+1} \stackrel{n_j}{\to} I_j)$ is defined in \eqref{equivO}. One has 
$$
\cO(I_{j+1} \stackrel{n_j}{\to} I_j)(f)(\lambda)=f(n_j\lambda)\qqq \lambda \in I_{j+1}
$$
and thus, using $n_j h_{j+1}=h_j$, one has for any $f\in \cO(I_j)$
$$
\rho(j+1,\cO(I_{j+1} \stackrel{n_j}{\to} I_j)(f))(\lambda)=\cO(I_{j+1} \stackrel{n_j}{\to} I_j)(f)(\lambda h_{j+1})=f(n_j\lambda h_{j+1})=f(\lambda h_j)=\rho(j,f)(\lambda).
$$
Thus the map $(j,f)\mapsto \rho(j,f)$ is compatible with the colimit and determines a map $\rho:\cO_{\ffp_H}\to \cR_H$
which is easily shown to be an isomorphism of semirings.\newline
 $(ii)$~Let  $(n_j)$, with $I_j$ be as in Lemma \ref{limlim2} and such that $H=\cup (\prod_1^j n_\ell)^{-1}\Z$.  We let, as above, $m_j:=\prod_1^{j-1} n_\ell$ so that  $H:=\cup\, m_j^{-1}\Z$.  We set $\iota(j,k):=k/m_j\in H$
so that  \eqref{iotak} follows by construction.
Then by Proposition \ref{proplimbis}, 
the stalk of the structure sheaf $\cO$ at the point $\ffq_H$ is the colimit $\cO_{\ffq_H}=\varinjlim \cO(I_j)$. We define a map $\delta:\cO_{\ffq_H}\to H$ as follows. We associate to $(j,f)$, with $f\in \cO(I_j)$, the element 
$\delta(j,f):=\iota(j,k)$ where $k=f'(0)\in \Z$ is the derivative of $f$ at $0\in I_j$. One then has 
$$
\delta(j+1,\cO(I_{j+1} \stackrel{n_j}{\to} I_j)(f))=\iota(j+1,(f(n_j\lambda))_{\lambda=0}'=\iota(j+1,n_jf'(0))=\iota(j,f'(0))=\delta(j,f).
$$
This shows that the map $\delta$ is well defined. Similarly, the equality $\alpha(j,f):=f(0)$ defines a map $\alpha:\cO_{\ffq_H}\to \rma$ and the pair $\rho=(\alpha,\delta)$ gives a map $\cO_{\ffq_H}\to \cZ_H$ which is both injective and surjective. One easily checks that this map is an isomorphism of $\cO_{\ffq_H}$ for the semiring structure whose multiplication corresponds to $(x,h)\bullet (x',h')=(x+x',h+h')$ and the addition to
$$
(x,h)\vee (x',h'):=\begin{cases} (x,h)\ \text{if}\ x>x'\\(x',h')\ \text{if}\ x'>x\\
(x,h \vee h') \ \text{if}\ x=x'.\end{cases}
$$ 
 \endproof
 Next, we describe the semiring $\cR_H$ of germs at $\lambda=1$ of $\rma$ valued piecewise affine (continuous)  convex functions $f(\lambda)$ with slopes in $H$. The germ of $f$ is  determined by the triple $(x,h_+,h_-)$, with $x\in \R$ and  $h_\pm\in H$ given by  $x=f(1)$, $f(1\pm \epsilon)=x
\pm h_\pm \epsilon$, for $\epsilon \geq 0$ small enough. The triples $(x,h_+,h_-)$ obtained from elements of $\cR_H$ are characterized by the condition $h_+\geq h_-$ which corresponds to the convexity of the function $f(1\pm \epsilon)=x
\pm h_\pm \epsilon$ for $\epsilon \geq 0$ small enough. The only other element of the semiring $\cR_H$ corresponds to the germ of the constant function $-\infty$. This  function plays the role of the ``zero" element for the following algebraic rules applied to the non-zero elements of the semiring. The ``addition" $\vee$ is given by the max between two germs and hence it is described by the formula
\begin{equation*}\label{additRh}
	(x,h_+,h_-)\vee (x',h'_+,h'_-):=\begin{cases} (x,h_+,h_-)\ \text{if}\ x>x'\\(x',h'_+,h'_-)\ \text{if}\ x'>x\\
(x,h_+\vee h'_+, h_-\wedge h'_-) \ \text{if}\ x=x'.\end{cases}
\end{equation*}
The ``product" of two germs is given by their sum and hence it is described by the formula
\begin{equation*}\label{prodRh}
(x,h_+,h_-)\bullet (x',h'_+,h'_-):=(x+x',h_++h'_+,h_-+h'_-).
\end{equation*}
The conditions on the functions $f(H)$ on subgroups $H\subset \R$ which define the structure sheaf are expressed locally in terms of the scaling flow by the two derivatives $D_\pm$ where 
\begin{equation}\label{dpm}
D_\pm(f)(H):=\lim_{\epsilon\to 0\pm}\frac{f((1+\epsilon)H)-f(H)}{\epsilon}	
\end{equation}
and can be written in the form
\begin{equation}\label{dpm1}
D_\pm(f)(H)\in H \qqq H	, \  \  D_+(f)(H)\geq D_-(f)(H) \qqq H.
\end{equation}
Indeed, taking a neighborhood of $H$ of the form $\{\lambda H\mid \lambda\in V\}$ where $V$ is a neighborhood of $1$, and considering the function $g(\lambda):=f(\lambda H)$ condition \eqref{dpm1} becomes 
$
\lambda\partial^\pm_\lambda g(\lambda)\in \lambda H
$, where the $\partial^\pm_\lambda$ are the directional derivatives. Thus the condition on the functions  becomes $\partial^\pm_\lambda g(\lambda)\in  H$ and
$\partial^+_\lambda g(\lambda)\geq\partial^-_\lambda g(\lambda)$ which is the characterization of the germs in a neighborhood of $H$. 

\subsection{The  points of $\scal2$ over $\rmax$}
The next  Theorem states that the process of extension of scalars from $\aarith$ to $\scal2$ does not affect the points over $\rmax$.

\begin{thm} \label{structure3} The canonical projection from the set $\scal2(\rma)$ of the points of the scaling site $\scal2$ defined over $\rmax$ 
to the points of the   topos $\rnt$ is bijective.
\end{thm}
The proof follows from the next Lemma
\begin{lem}\label{homtormax} $(i)$~The map $(x,h_+,h_-)\mapsto x$ is the only element of $\Hom_{\rma}(\cR_H,\rma)$.\newline
$(ii)$~The map $(x,h_+)\mapsto x$ is the only element of $\Hom_{\rma}(\cZ_H,\rma)$.
\end{lem}
\proof $(i)$~Let $\phi\in\Hom_{\rma}(\cR_H,\rma)$. First, we notice that the $\rma$-linearity shows that the image by $\phi$ of the constant function is $\phi(x,0,0)=x$. Then, we see that for any germ $f$ which is not identical to $-\infty$, there exists a constant function $g<f$. One then has 
$$
f\vee g=f\Rightarrow \phi(f)=\phi(f\vee g)=\phi(f)\vee \phi(g)\Rightarrow \phi(f)\geq \phi(g).
$$
This shows that one cannot have $\phi(f)=-\infty$. This argument shows that for any elements $f,g$ of $\cR_H$, one has  $f<g\Rightarrow \phi(f)\leq \phi(g)$ and it follows that
$$
x<x'\Rightarrow \phi(x,h_+,h_-)\leq\phi(x',h'_+,h'_-).
$$
Then,  since $\phi(x,0,0)=x$ one gets $\phi(x,h_+,h_-)=x$.
 
 $(ii)$~The proof is similar as for $(i)$.  \endproof

\subsection{The  sheaf of fractions and Cartier divisors}
Cartier divisors on a scheme  are defined as the global sections of the sheaf $\cK^\times/\cO^\times$ quotient of the sheaf of multiplicative groups of the  rings of fractions $\cK$  of the scheme  by the sub-sheaf $\cO^\times$ of invertible elements of the structure sheaf (\cf \cite{Hart} pp. 140-141). We  adapt this notion in this context of characteristic $1$.
\begin{prop}\label{propfunct} Let $\ffp_H$ be the point of the  topos $\rnt$ associated to a rank one subgroup $H\subset \R$.
The semiring $\cR_H$ of germs, at $\lambda=1$, of $\rma$-valued, piecewise affine, continuous  convex functions $f(\lambda)$ with slopes in $H$ is \mc and its semifield of fractions ${\rm Frac}\,\cR_H$ is the semifield of germs, at $\lambda=1$, of $\rma$-valued, piecewise affine, continuous functions $f(\lambda)$ with slopes in $H$ and endowed with the operations of the max and the addition of germs.
\end{prop}
\proof  To show that $\cR_H$ is \mc one considers an equation of the form $f+h=g+h$ for $f,g,h\in \cR_H$ and notice that since these functions take finite values at every point we get $f=g$. It follows that two pairs $(f,g)$ and $(h,k)$ of elements of $\cR_H$ define the same element of the associated semifield of fractions ${\rm Frac}\,\cR_H$ if and only if one has $f+k=g+h$. This is equivalent to write $f-g=h-k$ and hence to the equality of the  functions obtained as pointwise differences. Next, the addition in the semiring $\cR_H$ is given by the pointwise supremum and one needs to check that its extension to the associated semifield of fractions ${\rm Frac}\,\cR_H$ coincides with the pointwise supremum for  functions. This follows from the equality, valid for real numbers $x,y,z,t$ and obtained from translation invariance of $\vee$
$$
(x-y)\vee (z-t)=\left( (x+t)\vee (y+z)\right)-(y+t).
$$ \endproof 

As in the case of $\cR_H$ the  germ of a function $f$ is  characterized by the triple $(x,h_+,h_-)$ with $x\in \R$ and $h_\pm\in H$ given by  $x=f(1)$, $f(1\pm \epsilon)=x
\pm h_\pm \epsilon$, for $\epsilon \geq 0$ small enough. The triples $(x,h_+,h_-)$ obtained from elements of ${\rm Frac}\,\cR_H$ correspond to arbitrary values of $(x,h_+,h_-)\in \R\times H\times H$. As above the only other element of ${\rm Frac}\,\cR_H$ corresponds to the germ of the constant function $-\infty$ which plays the role of the ``zero" element. The algebraic rules for the other elements are given for the ``addition" $\vee$ of two germs  by the max of the two germs and hence by
\begin{equation}\label{additRhbis}
	(x,h_+,h_-)\vee (x',h'_+,h'_-):=\begin{cases} (x,h_+,h_-)\ \text{if}\ x>x'\\(x',h'_+,h'_-)\ \text{if}\ x'>x\\
(x,h_+\vee h'_+, h_-\wedge h'_-) \ \text{if}\ x=x'\end{cases}
\end{equation}
The ``product" of the germs is given by the sum of the two germs and hence by
\begin{equation}\label{prodRhbis}
(x,h_+,h_-)\bullet (x',h'_+,h'_-):=(x+x',h_++h'_+,h_-+h'_-).
\end{equation}

\subsection{The order at a point}
\begin{defn}\label{site2} Let $\ffp_H$ be the point of the $\rnt$ associated to a rank one subgroup $H\subset \R$ and let $f$ be an element in the stalk of $\cK$ at $\ffp_H$. Then,  the order of $f$ at $H$ is defined as $\ord(f) = h_+-h_-\in H\subset \R$, where $h_\pm = \lim_{\epsilon\to 0\pm}\frac{f((1+\epsilon)H)-f(H)}{\epsilon}$.
\end{defn}

\begin{prop}\label{ordercomp} For any two elements $f,g$ in the stalk of $\cK$ at $\ffp_H$ one has
  \begin{equation}\label{ordcomp1}
\ord(f\vee g)\geq \ord(f)\wedge \ord(g)
\end{equation}
\begin{equation}\label{ordcomp2}
\ord(f+ g)= \ord(f)+ \ord(g).
\end{equation}
\end{prop}
\proof The germ of $f\vee g$ is given by the max of the two germs and hence by
\eqref{additRhbis}. 
When $x\neq x'$ the  inequality \eqref{ordcomp1} follows from $h\geq h\wedge h'$, $h'\geq h\wedge h'$. When $x=x'$ one has 
$\ord(f\vee g)=(h_+\vee h'_+)-(h_-\wedge h'_-)$. Thus the inequality \eqref{ordcomp1} follows from the general fact
$$
(a\vee b)-(c\wedge d)\geq (a-c)\vee(b-d)\geq (a-c)\wedge (b-d)\qqq a,b,c,d\in \R.
$$
Note that one needs the $\wedge$ in \eqref{ordcomp1} to take care of the cases $x\neq x'$, but when $x=x'$ the $\vee$  works. 

Finally, the germ of $f +
g$ is given by the sum of the two germs and hence \eqref{ordcomp2} holds. \endproof

\section{The periodic orbits $C_p$}

Let $p$ be a prime and consider the subspace $C_p$ of points of  $\rnt$ corresponding to subgroups $H\subset \R$ which are abstractly isomorphic to the subgroup $H_p\subset \Q$ of fractions with denominator a power of $p$. In this section we study the (quasi)-tropical structure of these curves, we develop the theory of theta-functions and finally we formulate a Riemann-Roch problem and establish a Riemann-Roch formula.

\begin{lem}\label{periodp} $(i)$~The map $\R_+^*\to C_p$, $\lambda\mapsto \lambda H_p$  induces the isomorphism $\eta_p: \R_+^*/p^\Z\to C_p$.\newline
$(ii)$~The pull-back by $\eta_p$ of the (restriction to $C_p$ of the) structure sheaf $\cO$ is the sheaf $\cO_p$ on $\R_+^*/p^\Z$ of piecewise affine (continuous), convex functions with slopes in $H_p$.
\end{lem}
\proof $(i)$~For $\lambda\in \R_+^*$ one has $\lambda H_p\in C_p$. An abstract isomorphism of ordered groups $\phi:H_p\to H\subset \R$ is given by multiplication by $\phi(1)=\lambda$. Since $pH_p=H_p$ the  map $\R_+^*\to C_p$ induces a surjective map $\eta_p: \R_+^*/p^\Z\to C_p$. We show that $\eta_p$ is injective. If $\lambda H_p=\lambda' H_p$, then $\mu H_p=H_p$ for $\mu=\lambda/\lambda'$. Thus $\mu=a/p^n\in H_p$ and the same result holds for $\mu^{-1}$, thus $\mu$ is a power of $p$.

$(ii)$~The condition on the functions $f(H)$ on subgroups $H\subset \R$ which defines the structure sheaf are expressed locally in terms of the scaling flow by the two derivatives $D_\pm$ of \eqref{dpm} and can be written in the form \eqref{dpm1}. On the function $g=f\circ \eta_p$,  
 $g(\lambda):=f(\lambda H_p)$, condition \eqref{dpm1} becomes 
$
\lambda\partial^\pm_\lambda g(\lambda)\in \lambda H_p
$, $\lambda\partial^+_\lambda g(\lambda)\geq\lambda\partial^-_\lambda g(\lambda)$ (where the $\partial^\pm_\lambda$ are the directional derivatives) or equivalently: $\partial^\pm_\lambda g(\lambda)\in  H_p$ and $\partial^+_\lambda g(\lambda)\geq\partial^-_\lambda g(\lambda)$ thus one gets $(ii)$.\endproof

\subsection{Divisors}
\begin{prop}\label{cartierdiv} $(i)$~The sheaf of quotients of the sheaf of semirings $\cO_p$ is the sheaf $\cK_p$ on $\R_+^*/p^\Z$ of piecewise affine, continuous  functions with slopes in $H_p$,  endowed with the two operations max, $+$. \newline
$(ii)$~The quotient sheaf of Cartier divisors $\cdiv(C_p):=\cK_p^\times/\cO_p^\times$ is isomorphic to the sheaf $\mathscr{D}iv(C_p)$ of naive divisors, \ie of maps  $H\mapsto D(H)\in H$, such  that
$$
\forall \lambda\in\R_+^*, \ \exists V \ \text{open} \ \lambda\in V: \  D(\mu)=0\qqq \mu \in V, \ \mu\neq \lambda.
$$
\end{prop}
\proof $(i)$~The proof is the same as for Proposition \ref{propfunct}.

$(ii)$~A germ $(x,h_+,h_-)\in \cR_H$ is invertible for the multiplicative structure given by \eqref{prodRhbis} iff one has $h_+=h_-$, thus the map given by the order as in Definition \ref{site2} is an isomorphism  \endproof

\begin{defn}\label{defndiv} A {\em divisor} on $C_p$ is a global section of the sheaf 
$\cdiv(C_p)\simeq \mathscr{D}iv(C_p)$. The collection of these global sections is denoted by $\div(C_p)$.
\end{defn} 
By Proposition \ref{cartierdiv}, a divisor on $C_p$ is uniquely specified by a global section of $\mathscr{D}iv(C_p)$  and by compactness, these are the maps $H\in C_p\mapsto D(H)\in H$ with finite support, where the support is defined 
	\begin{equation*}\label{support}
	\supp(D):=\{H\mid D(H)\neq 0\}.
	\end{equation*}
	In other words a divisor $D$ on $C_p$ is a section, vanishing everywhere except on a finite subset of $C_p$, of the projection on the base from the total space of the bundle formed of pairs $(H,h)$ where $H\subset \R$ is a subgroup abstractly isomorphic to the subgroup $H_p\subset \Q$ and $h\in H$.
	The sheaf $\cK_p$ has global sections and they form the semifield  $\cK(C_p):= H^0(\R_+^*/p^\Z,\cK_p)$.  
	
	\begin{prop}\label{propdiv}
	$(i)$~The divisors $\div(C_p)$ form an abelian group under pointwise addition.\newline
	$(ii)$~The condition $D'(H)\geq D(H)$, $\forall H\in C_p$, defines a partial order on the group $\div(C_p)$.\newline
	$(iii)$~The following map defines a surjective group homomorphism compatible with the partial order 
	\begin{equation*}\label{degg}
\deg: \div(C_p)\to \R,\quad \deg(D):=\sum_H D(H)\in\R.
\end{equation*}
	$(iv)$~The map which associates to  $f\in \cK^\times(C_p)$  the principal divisor 
	\begin{equation}\label{princdiv}
	(f):=\sum_H (H,\ord_H(f))
	\end{equation} 
	determines a group homomorphism $ \cK^\times(C_p)\to \div(C_p)$.\newline
	$(v)$~The subgroup $\cP\subset \div(C_p)$ of principal divisors is contained in the 
	kernel of $\deg:\div(C_p)\to \R$ 	
	\begin{equation*}\label{order}
\sum_{H} \ord_H(f)=0\qqq f\in \cK(C_p).
\end{equation*}
	\end{prop}
	\begin{figure}
\begin{center}
\includegraphics[scale=0.7]{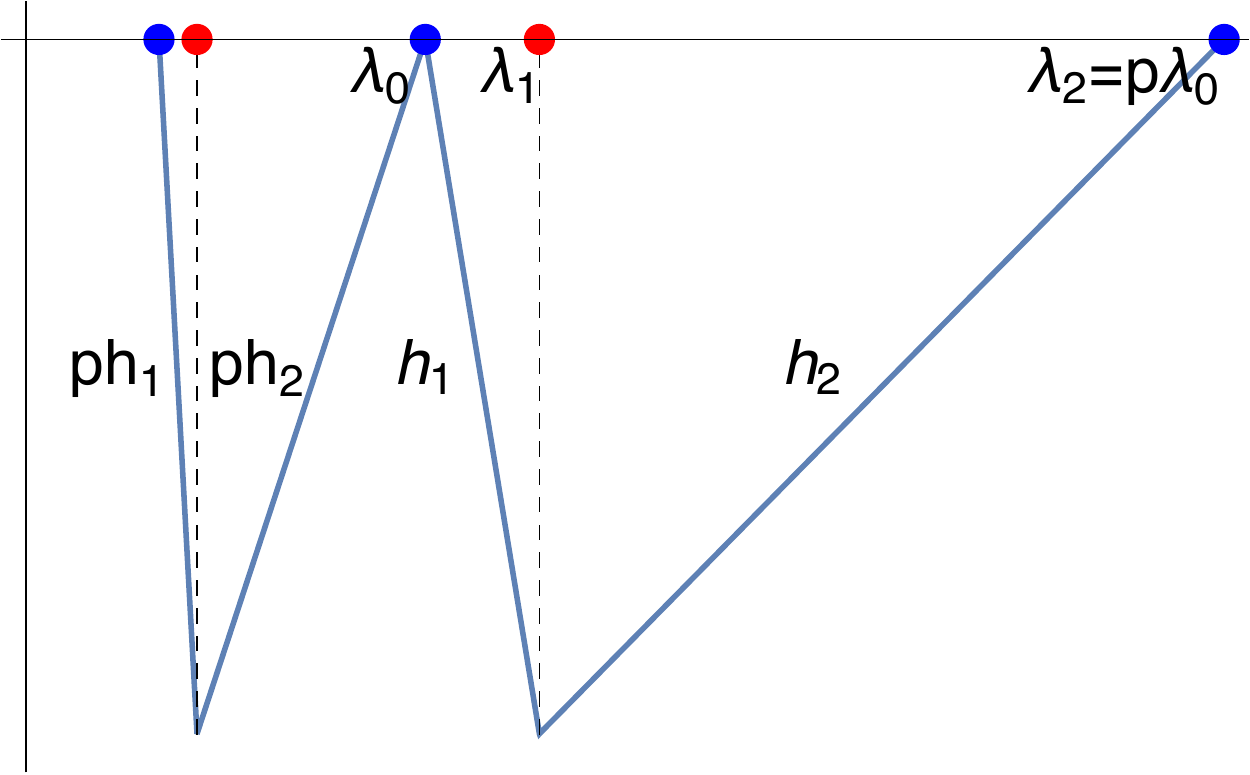}
\end{center}
\caption{A function $f\in\cK^\times(C_p)$. \label{tropfct} }
\end{figure}
	\proof $(i)$~By construction $\div(C_p)$ is the direct sum of the groups $H$, with $H\in C_p$.\newline
	$(ii)$~Each group $H$ is  ordered  and the condition $D'(H)\geq D(H)$, $\forall H\in C_p$ defines the natural partial order on their direct sum.\newline
	$(iii)$~The statement follows from the fact that the groups $H$ are subgroups of $\R$ and their union  is $\R$. \newline
	$(iv)$~For any $f\in \cK^\times(C_p)$ the sum in \eqref{princdiv} is finite and defines a divisor. Moreover the formula \eqref{ordcomp2} shows that one obtains a group homomorphism.\newline
	 $(v)$~A global section $f\in \cK(C_p)$ is a real valued function $f:\R_+^*\to \R$ which is  piecewise affine, continuous,   with slopes in $H_p$ and fulfills $f(p\lambda)=f(\lambda)$ $\forall \lambda\in\R_+^*$ (\cf~Figure~\ref{tropfct}). Such a function is uniquely determined by its restriction to the fundamental domain $[\lambda_0,p\lambda_0]$ and the only  constraint on this restriction is the periodicity: $f(p\lambda_0)=f(\lambda_0)$. Let $\lambda_0<\lambda_1<\ldots <\lambda_{n-1}<\lambda_n=p\lambda_0$ be a finite sequence of positive real numbers such that $f$ is affine with slope $h_j\in H_p$ on the interval $I_j=[\lambda_{j-1},\lambda_j]$. One has for any $j\in \{1,\ldots ,n\}$
$$
f(\lambda_j)=f(\lambda_0)+\sum_1^j (\lambda_i-\lambda_{i-1})h_i.
$$
By applying the equality above when $j=n$, and using the periodicity property $f(\lambda_n)=f(p\lambda_0)=f(\lambda_0)$ one obtains
$$
\sum_1^n (\lambda_i-\lambda_{i-1})h_i=0.
$$
By using $\lambda_n=p\lambda_0$ one also has
$$
(\lambda_1-\lambda_0)h_1+(\lambda_2-\lambda_1)h_2+\ldots +(\lambda_n-\lambda_{n-1})h_n=\sum_1^{n-1}  \lambda_i(h_i-h_{i+1})+\lambda_0(ph_n-h_1).
$$
For $1\leq i\leq n-1$ one has $\lambda_i(h_i-h_{i+1})=-{\rm Order}(f)(\lambda_i)$ where we set ${\rm Order}(f)(\lambda):=\ord_{\lambda H_p}f$. Moreover, notice that $ph_n$ is the slope of the function $f$ in the interval $(1/p) I_n$ so that the order of $f$ at $\lambda_0$ is $\lambda_0(h_1-ph_n)$. Then, the above equality shows that the sum of all orders must vanish. \endproof

Next, we consider the problem of constructing a global section $f\in \cK(C_p)$ whose divisor is an assigned divisor of degree zero. 

Let $J_p:=(p-1)H_p \subset H_p$ be the principal ideal generated by the integer $(p-1)$ in the ring $H_p$. One has the following exact sequence of rings 
\begin{equation*}\label{rings}
	0\to J_p\to H_p\stackrel{\chi}{\to} \Z/(p-1)\Z\to 0
\end{equation*}
where $\chi(a/p^n):=a$ mod. $(p-1)$, for any $a\in \Z$ and $n\in \N$. 
\begin{prop}\label{propdiv1}
	$(i)$~Let $H\in C_p$, then there exists a unique map $\chi_H:H\to H_p/(p-1)H_p\simeq \Z/(p-1)\Z$ such that for any  $\lambda\in \R_+^*$ with $H=\lambda H_p$, one has $\chi_H=\chi\circ \lambda^{-1}$.\newline
	$(ii)$~The group $\cP$ of principal divisors  is contained in the kernel of the group homomorphism 
\begin{equation*}\label{jacmap}
\chi: \div(C_p)\to Z/(p-1)\Z,\quad \chi(D):=\sum_H \chi_H(D(H)).
\end{equation*}
\end{prop}
\proof $(i)$~Given $H\in C_p$ the $\lambda 's \in \R_+^*$ such that $H=\lambda H_p$ provide maps $\lambda^{-1}:H\to H_p$ which differ from each other by multiplication by a power of $p$, thus the corresponding map $H\to H_p/(p-1)H_p\simeq \Z/(p-1)\Z$ is independent of the choice of $\lambda$. \newline
$(ii)$~We use the notations of the proof of Proposition \ref{propdiv} $(v)$. The support of the principal divisor $(f)$ is contained in $\{\lambda_j H_p\mid 0\leq j\leq n-1\}$. For $1\leq j\leq n-1$, one has ${\rm Order}(f)(\lambda_j)=\lambda_j(h_{j+1}-h_j)$ and thus 
$$
\chi_{\lambda_jH_p}(\ord_{\lambda_jH_p}(f))=\chi(h_{j+1}-h_j)=\chi(h_{j+1})-\chi(h_j).
$$ 
For $j=0$: ${\rm Order}(f)(\lambda_0)=\lambda_0(h_{1}-ph_n)$ and thus $$\chi_{\lambda_0H_p}(\ord_{\lambda_0H_p}(f))=\chi(h_{1}-ph_n)=\chi(h_1)-\chi(h_n).
$$
Thus one gets $\chi((f))=0$ as required. \endproof 

\begin{thm}\label{thmjaccp} The map defined by the pair of homomorphisms  
\begin{equation}\label{jacmapbis}
(\deg,\chi):\div(C_p)/\cP\to \R\times (\Z/(p-1)\Z)
\end{equation}
is an isomorphism of abelian groups. 
\end{thm}
\proof We show first that any divisor $D$ such that $\deg(D)=0$ and $\chi(D)=0$ is principal.   Indeed, the divisor $D$ can be written in the form
$$
D=\sum_{0\leq j\leq n-1}(\lambda_j H_p, a_j\lambda_j), \qquad \lambda_0<\lambda_1<\ldots <\lambda_{n-1}<p\lambda_0
$$
where the elements $a_j\in H_p$ are such that $\sum_0^{n-1}a_j\lambda_j=0$. Thus one is given the intervals $I_j=[\lambda_{j-1},\lambda_j]$ with $\lambda_n=p\lambda_0$. In order to show that $D$ is principal one needs to find elements $h_j\in H_p$ for $1\leq j\leq n$, giving the slope of $f$ on $I_j$ and hence such that 
$$
a_0=h_1-ph_n,\ a_1=h_2-h_1, \ a_2=h_3-h_2, \ \ldots ,a_{n-1}=h_n-h_{n-1}.
$$
The solution of the above system of equations is unique and  of the form 
$$
h_n=\sigma/(1-p),\ h_{n-1}=h_n-a_{n-1}, \ldots ,h_1=h_2-a_1, \quad \sigma :=\sum a_j.
$$
Thus it takes values in $H_p$ if and only if the element $\sigma$ is divisible by
$p-1$ in $H_p$ \ie iff $\chi(D)=0$. If this is the case one defines the function $f$ to be affine with slope $h_j$ in the interval $I_j=[\lambda_{j-1},\lambda_j]$ and normalized by $f(\lambda_0)=0$.
One has, using $\lambda_n=p\lambda_0$
$$
f(\lambda_n)=\sum_1^n(\lambda_j-\lambda_{j-1})h_j=-\lambda_0 h_1+\sum_1^{n-1}\lambda_j(h_j-h_{j+1})+\lambda_n h_n=-\sum_0^{n-1}a_j\lambda_j=0.
$$
This argument shows that one can extend $f$ by periodicity so that $f(p\lambda)=f(\lambda)$ and obtain a continuous, piecewise affine function with slopes in $H_p$. Moreover, by construction the divisor of $f$ is $D$, since the discontinuities of the derivative take place at the $\lambda_j$'s and are given by the $a_j$'s. 
 It remains to show that the restriction of the map $\chi$ to the subgroup $\div_0(C_p)\subset \div(D)$ of divisors of degree $0$ is surjective  onto $\Z/p\Z$.  Given $m\in\Z/(p-1)\Z$, it is enough to find an increasing finite sequence of real numbers $\lambda_j>0$, $\lambda_0<\lambda_1<\ldots<\lambda_n=p\lambda_0$ and elements $a_j\in H_p$ such that $\sum_0^{n-1}a_j\lambda_j=0$ and  $\sum_0^{n-1}a_j=m$ mod. $(p-1)H_p$. Both the group $H_p\subset \R$ and the subset $H_p^{(m)}:=\{x\in H_p\mid x=m \ \text{mod.}\ (p-1)H_p\}$ are dense in $\R$. This shows that one can fix arbitrarily the $(\lambda_j)$ and all the $a_j$ except for one of them, say $a_k$, and then choose $a_k\in H_p$ with preassigned $\chi(a_k)\in \Z/(p-1)\Z$ such that for a choice of $\lambda'_k$ close to $\lambda_k$ one gets 
$$
\sum_0^{n-1}a_j\lambda'_j=0, \qquad  \sum_0^{n-1}\chi(a_j)=m \in \Z/(p-1)\Z.
$$
The same argument also shows that one can arbitrarily prescribe both $\deg(D)$ and $\chi(D)$ for divisors $D\in \div(C_p)$.
Since the group law on divisors is  given by pointwise addition of sections, both maps $\deg:\div(C_p)\to \R$ and $\chi:\div(C_p)\to \Z/(p-1)\Z$ are group homomorphisms and one obtains the isomorphism of groups \eqref{jacmapbis}.\endproof

\subsection{Symmetries}\label{sectsymm}
The scaling site $\scal2$ inherits from its construction by extension of scalars, as in Section \ref{sectextS}, several symmetries such as the arithmetic Frobenius associated to the automorphisms of $\rmax\sim \rma$ or the absolute Frobenius associated to the frobenius endomorphisms of the semirings of characteristic one of the structure sheaf. These symmetries induce 
 symmetries of the curves $C_p$ and in this subsection we describe them as operators acting on the structure sheaf $\cO_p$ of these curves. Since the divisors on $C_p$ are  Cartier divisors, one then obtains, as a byproduct,  induced actions on 
$\div(C_p)$ and on its quotient $\div(C_p)/\cP$.

\subsubsection*{Arithmetic Frobenius}\label{sectaritfrob}
 The analogue of the arithmetic Frobenius $\fr^a$ on sections $f$ of $\cO_p$ is defined as follows 
\begin{equation*}\label{arithfrob}
\fr^a_\mu(f)(\lambda):=\mu \, f(\mu^{-1}\lambda)\qqq \lambda, \mu\in \R_+^*	
\end{equation*}
This operator preserves the properties of $f$ ($f$ is  convex, piecewise affine with slopes in $H_p$ and periodic $f(p\lambda)=f(\lambda)$) as well as the algebraic operations $(\vee,+)$. Thus it defines automorphisms denoted by $\fr^a_\mu$. On a function $f$ which is locally written as $\vee (n_j\lambda +a_j)$ the action of $\fr^a_\mu$ replaces the $a_j$'s by $\mu a_j$ and this corresponds to the arithmetic Frobenius. The induced action on the stalks of $\cO_p$ associates to a germ at $\lambda$ a germ at $\mu\lambda$ given by $x\mapsto \mu \, f(\mu^{-1}x)$, for $x$ near $\mu\lambda$. One thus obtains a morphism from germs at $H$ to germs at $\mu H$ 
\begin{equation}\label{arithfrob1}
\fr^a_\mu:\cR_H\to \cR_{\mu H}, \qquad  (x,h_+,h_-)\mapsto (\mu x,\mu h_+,\mu h_-).
\end{equation}
Indeed, one has $f((1+\epsilon)H)\sim f(H)+\epsilon h_+$, for $\epsilon>0$ small and thus 
$$
\fr^a_\mu(f)((1+\epsilon)\mu H)=\mu \, f((1+\epsilon)H)\sim \mu f(H)+\mu \epsilon h_+.
$$
Next, we list some properties of the induced action by $\fr^a_\mu$ on divisors. We use the notation 
$D=\sum (H_j,h_j)$, with $h_j\in H_j$ for the divisor with support on the set $\{H_j\}$ and such that $D(H_j)=h_j$ $\forall j$. 

\begin{lem}\label{arithfrob2} The  action induced by $\fr^a_\mu$ on divisors is given by 
$$
D=\sum (H_j,h_j)\mapsto \fr^a_\mu(D)=\sum (\mu H_j,\mu h_j).
$$
	This action preserves the homomorphism $\chi:\div(C_p)\to \Z/(p-1)\Z$ as well as the subgroup $\cP$ of principal divisors and it acts on the degree of a divisor by multiplication by $\mu$.
\end{lem}
\proof The first part of the statement follows from \eqref{arithfrob1} which shows that the singularity of the derivative of $f$ at $H$ gives a singularity of the derivative of $\fr^a_\mu(f)$ at $\mu H$, while the order is multiplied by $\mu$. The value of $\chi(H_j,h_j)$ is obtained as $\chi(k)$ where $H_j=\lambda H_p$, $h_j=\lambda k$. One can choose the same $k$ for the pair $(\mu H_j,\mu h_j)$ and this implies that $\chi(\fr^a_\mu(D))=\chi(D)$. The degree of $D=\sum (H_j,h_j)$ is $\deg(D)=\sum h_j$ and the degree of $\fr^a_\mu(D)$ is $\sum \mu h_j=\mu \deg(D)$. Since $\cP=\Ker(\deg,\chi)$ this group is preserved by the action of $\fr^a_\mu$.\endproof

\subsubsection*{Relative Frobenius}\label{sectrelfrob}

 Let $h\in H_p$, $h>0$.  One defines the operator $\fr^r_h$ acting on functions  by
\begin{equation*}\label{relatfrob}
\fr^r_h(f)(\lambda):= f(h\lambda)\qqq \lambda\qqq h\in H_p^+.	
\end{equation*}
This operator acts on a section $f$ of $\cO_p$ locally written as $f=\vee (n_j\lambda +a_j)$ by replacing the $n_j$ by $hn_j$, while leaving the $a_j$ unchanged. In particular it is $\rma$-linear.
Since the powers $p^n$ act trivially in view of the periodicity of $f$, the operator $\fr^r_h$ only depends upon the class of $h$ in the quotient $H_p^+/p^\Z$ which is the multiplicative monoid $\ntp$ of positive integers relatively prime to $p$. The induced action on the stalks associates to a germ at $\lambda$,  the germ at $h^{-1}\lambda$ given by $x\mapsto  f(hx)$ for $x$ near $h^{-1}\lambda$. One thus obtains a morphism from germs at $H$ to germs at $h^{-1} H$ 
\begin{equation}\label{relatfrob1}
\fr^r_h:\cR_H\to \cR_{h^{-1} H}, \qquad  (x,h_+,h_-)\mapsto ( x,h_+, h_-)\in \rma\times h^{-1} H\times h^{-1} H
\end{equation}
as follows from the identity $\fr^r_h(f)((1+\epsilon)h^{-1} H)=f((1+\epsilon)H)$.
\begin{lem}\label{relatfrob2} The  action induced by $\fr^r_h$ on divisors is given by 
$$
D=\sum (H_j,h_j)\mapsto \fr^r_h(D)=\sum (h^{-1}H_j, h_j).
$$
	This action preserves the degree $\deg\circ \fr^r_h=\deg$,  the subgroup $\cP$ of principal divisors and it acts on the invariant $\chi$ by multiplication by $\chi(h)\in\Z/(p-1)\Z$.
\end{lem}
\proof The first part of the statement follows from \eqref{relatfrob1}. One has $\deg(\fr^r_h(D))=\sum h_j=\deg(D)$. The value of $\chi(H_j,h_j)$ is obtained as $\chi(k)$ where $H_j=\lambda H_p$, $h_j=\lambda k$. One then has $h^{-1}H_j=h^{-1}\lambda H_p$ and $h_j=(h^{-1}\lambda)h k$,
so that  $\chi(h^{-1}H_j,h_j)=\chi(hk)=\chi(h)\chi(k)$. This argument also shows that $\cP=\Ker(\deg,\chi)$ is preserved by the action of $\fr^r_h$.\endproof

\subsubsection*{Absolute Frobenius}\label{sectabsfrob}

This operator acts on sections $f$ of $\cO_p$ by composition with the Frobenius of $\rma$ which, in this logarithmic notation, multiplies the function $f$ by a constant. In order to preserve the property of the slopes in $H_p$ one takes this constant to be an element in $H_p$. More precisely, the action of the absolute Frobenius on functions is given, for any $h\in H_p^+$, by the formula
\begin{equation*}\label{absfrob}
\fr_h(f)(\lambda):= h\,f(\lambda)\qqq \lambda\qqq h\in H_p^+.	
\end{equation*}
Its properties follow from the properties of the two previous operators since 
\begin{equation}\label{absfrob1}
\fr_h= \fr^a_h\circ \fr^r_h=\fr^r_h\circ \fr^a_h\qqq h\in H_p^+.	
\end{equation}
Indeed, one has $\fr^a_h\circ \fr^r_h(f)(\lambda)=hf(h^{-1}h\lambda)=hf(\lambda)$ and similarly $\fr^r_h\circ \fr^a_h=\fr_h$.
\begin{lem}\label{absfrob2} The  action induced by $\fr_h$ on divisors is given by 
$$
D=\sum (H_j,h_j)\mapsto \fr_h(D)=\sum (H_j, hh_j).
$$
	This action is trivial on the points, \ie it fixes $\supp(D)$. It preserves  the subgroup $\cP$ of principal divisors,  acts on the invariant $\chi$ by multiplication by $\chi(h)\in\Z/(p-1)\Z$ and on the degree by multiplication by $h\in H_p^+\subset \R_+^*$.
\end{lem}
\proof This follows from \eqref{absfrob1} using Lemmas \ref{arithfrob2} and \ref{relatfrob2}.
\endproof
Notice that each point $H\in C_p$ determines a morphism 
$$
f\in \cK(C_p)\mapsto p_H(f)=f(H)\in \rma\simeq \rmax
$$
and the absolute Frobenius is compatible with this morphism: $p_H(\fr_h(f))=\fr_h(p_H(f))$.
In particular it acts trivially on the spectrum.

\subsection{Theta functions}
In this subsection we provide an explicit construction of elements of $\cK(C_p)$  with assigned divisor, by introducing (tropical) analogues of theta functions. 
We proceed by analogy with the construction of theta functions for elliptic curves $E_t(k)=k^\times/t^\Z$ over non-archimedean local fields $k$ as in \cite{Tate}. In that case the theta function is defined as 
\begin{equation}\label{tatetheta}
\theta(w,t)=\sum_\Z (-1)^n t^{\frac{n^2-n}{2}}	w^n=(1-w)\prod_1^\infty(1-t^m)(1-t^m w)(1-t^m w^{-1})
\end{equation}
and satisfies the functional equation: $-w\theta(tw,t)=\theta(w,t)$. Its relation with the standard theta function $\vartheta _1(u,q)$  is given by the equation: $\vartheta _1(u,q)=i\sqrt[4]{q}\, e^{-iu}\theta(e^{2iu},q^2)$. 
This suggests to transpose \eqref{tatetheta} naively in order to define theta functions for $C_p=\R_+^*/p^\Z$. In our framework, the role of the product is replaced by addition, thus the following infinite sums replace the infinite products on the left
$$
\prod_0^\infty (1-t^m w)\rightsquigarrow f_+(\lambda):=\sum_0^\infty \left(0 \vee (1-p^{m}\lambda)\right)
$$
$$
\prod_1^\infty (1-t^m w^{-1})\rightsquigarrow f_-(\lambda):=\sum_1^\infty \left(0 \vee (p^{-m}\lambda-1)\right).$$

\begin{figure}
\begin{center}
\includegraphics[scale=0.7]{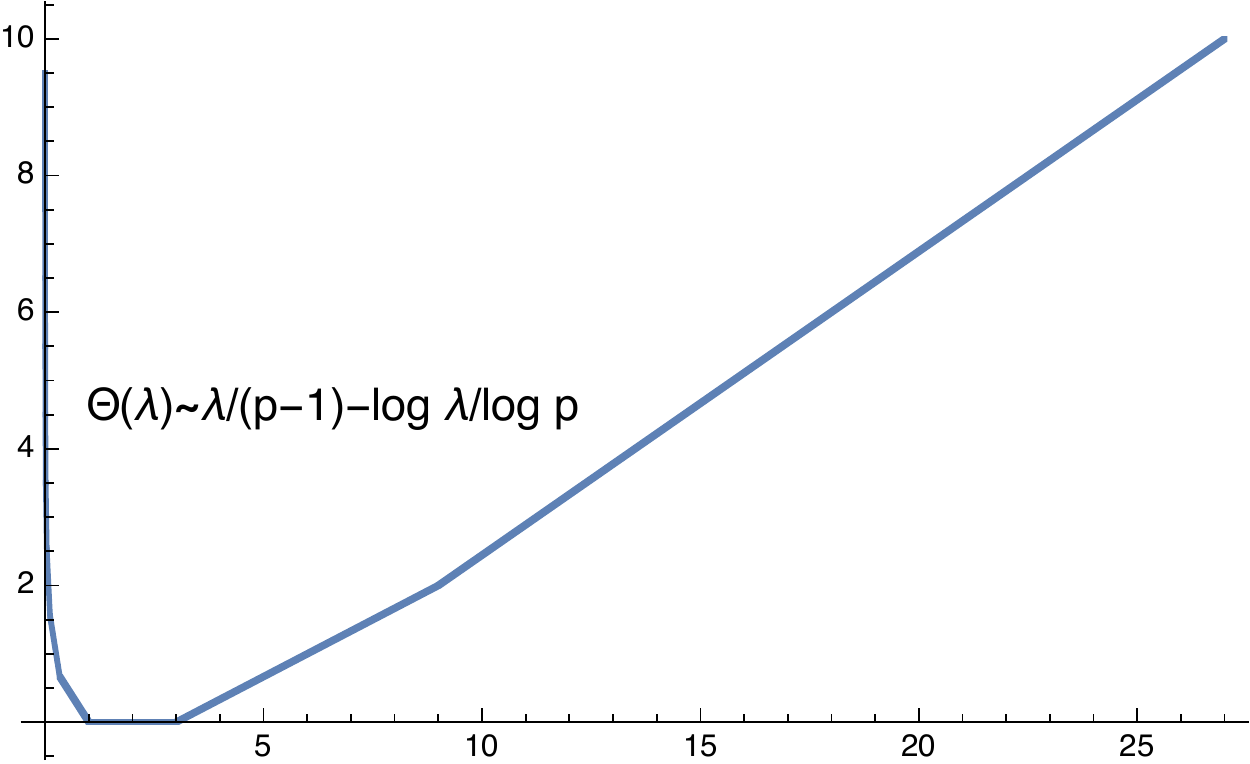}
\end{center}
\caption{The function $\theta$ for $p=3$. \label{thetafct} }
\end{figure}

\begin{lem}\label{lemtheta1} $(i)$~The functions $f_\pm(\lambda)$, $\lambda\in (0,\infty)$  are convex (continuous), piecewise affine with slopes in $H_p$ and 
\begin{equation}\label{theta1}
	f_+(p\lambda)-f_+(\lambda)=-\left(0 \vee (1-\lambda)\right), \qquad  
	f_-(p\lambda)-f_-(\lambda)=\left(0 \vee (\lambda-1)\right).
\end{equation}	
$(ii)$~The function $\theta(\lambda):=f_+(\lambda)+f_-(\lambda)$, $\lambda\in (0,\infty)$ (\cf~Figure~\ref{thetafct}), is convex (continuous), piecewise affine with slopes in $H_p$ and fulfills the equation: $\theta(p\lambda)=\theta(\lambda)+\lambda-1$, $\forall \lambda\in (0,\infty)$.\newline
$(iii)$~One has 
$$
\vert\theta(\lambda)-\left(\frac{1}{p-1}\lambda -\log \lambda/\log p\right)\vert \leq 1 \qquad \forall \lambda\in (0,\infty).
$$
\end{lem}
\proof $(i)$~Notice that the sum $\sum_0^\infty \left(0 \vee (1-p^{m}\lambda)\right)$ has only finitely many non-zero terms since $p^{m}\lambda>1$ for $m$ large enough. Each of these terms is
convex continuous, piecewise affine with slopes in $H_p$ and thus the same property holds for $f_+$. Moreover, the difference $f_+(p\lambda)-f_+(\lambda)$ is given by the single term corresponding to $m=0$,  \ie  $f_+(p\lambda)-f_+(\lambda)=-\left(0 \vee (1-\lambda)\right)$. Similarly for 
$f_-(\lambda)=\sum_1^\infty \left(0 \vee (p^{-m}\lambda-1)\right)$ one gets the term $\left(0 \vee (\lambda-1)\right)$ in $f_-(p\lambda)-f_-(\lambda)$. Thus  we obtain \eqref{theta1}.\newline 
$(ii)$~The first part of the statement follows from $(i)$.  Moreover one has 
$$
\theta(p\lambda)-\theta(\lambda)=f_+(p\lambda)-f_+(\lambda)+f_-(p\lambda)-f_-(\lambda)
=-\left(0 \vee (1-\lambda)\right)+\left(0 \vee (\lambda-1)\right)=\lambda-1
$$
since for any real number $x$ one has $(0\vee x)-(0\vee -x)=x$.\newline
$(iii)$~Let $g(\lambda):=\frac{1}{p-1}\lambda -\log \lambda/\log p$. Then 
$
g(p\lambda)-g(\lambda)=\lambda-1
$. 
Thus the function $k(\lambda)=\theta(\lambda)-g(\lambda)$ fulfills $k(p\lambda)=k(\lambda)$. This periodicity shows that $\vert k(\lambda)\vert\leq\max_{[1,p]}\vert k(u)\vert$. Since $\theta(u)=0$, $\forall u \in [1,p]$ one just needs to check that $\vert g(u)\vert \leq 1$,  $\forall u \in [1,p]$. In this interval the convex function $g$ varies between its value at the end points: $\frac{1}{p-1}$ and its minimum $g(\lambda)$ at $\lambda=\frac{p-1}{\log p}\in [1,p]$, whose value 
is $g(\frac{p-1}{\log p})=(1-\log(p-1)+\log\log p)/\log p\geq -1$.
\endproof

We define, for $h\in H_p, h>0$ and $\mu\in \R_+^*$, the function
\begin{equation*}\label{defnthet}
	\Theta_{h,\mu}(\lambda):=\mu\, \theta(\mu^{-1}h\lambda). 
\end{equation*}
It is a convex continuous, piecewise affine function with slopes in $H_p$ and fulfills the equation
\begin{equation}\label{propthet}
	\Theta_{h,\mu}(p\lambda)=\Theta_{h,\mu}(\lambda)+h\lambda -\mu 
\end{equation}
since by Lemma \ref{lemtheta1} $(ii)$ one has 
$$
\mu\, \theta(\mu^{-1}hp\lambda)=\mu\, \theta(\mu^{-1}h\lambda)+\mu\,\mu^{-1}h\lambda-\mu.
$$
It follows that
\begin{equation*}\label{propthet1}
	\Theta_{ph,\mu}(\lambda)=\Theta_{h,\mu}(\lambda)+h\lambda -\mu. 
\end{equation*}
To a pair of labels $(h,\mu)$ we associate the divisor which is everywhere zero except on the subgroup $H=\mu h^{-1} H_p$ (abstractly isomorphic to $H_p$) where it takes the value $\mu\in H$.  We let
\begin{equation*}\label{defndivthet}
	\delta(h,\mu):=(\mu h^{-1} H_p, \mu).
\end{equation*}
Note that by construction one has $\delta(ph,\mu)=\delta(h,\mu)$, so that $\delta(h,\mu)$ only depends upon the class of $h$ modulo the multiplication by powers of $p$. Such a class is uniquely specified as that of an integer $m>0$, prime to $p$. 
Given a simple divisor, \ie a positive divisor supported on a single $H$, and $\mu\in H$, $\mu>0$ one can find $\lambda>0$ such that $H=\lambda H_p$ and thus $h\in H_p$ with $\mu =h\lambda$. One then has  $\delta(h,\mu)=(H,\mu)$. The choice of $\lambda$ (and hence of $h$) is unique only up-to multiplication by a power of $p$. 

The classical description of elliptic functions in terms of theta functions admits the following counterpart in our framework, which provides, in particular, with another proof of Theorem \ref{thmjaccp} by an explicit construction of an $f\in \cK(C_p)$ with given divisor.
\begin{prop}\label{proptheta1}
Let $D=D_+-D_-\in \div(C_p)$ be a divisor ($D_\pm\in \div^+$)  and $(h_i,\mu_i)\in H_p^+\times \R_+^*$, $(h'_j,\mu'_j)\in H_p^+\times \R_+^*$ such that $D_+=\sum \delta(h_i,\mu_i)$ and $D_-=\sum \delta(h'_j,\mu'_j)$. Then, if $\deg(D)=0$ and $h\in H_p$ fulfills $(p-1)h=\sum h_i-\sum h'_j$, the following function  
\begin{equation}\label{propthet2}
	f(\lambda):=\sum_i \Theta_{h_i,\mu_i}(\lambda)-\sum_j \Theta_{h'_j,\mu'_j}(\lambda)-h\lambda 
\end{equation}	
is continuous, piecewise affine with slopes in $H_p$, fulfills $f(p\lambda)=f(\lambda)$ $\forall \lambda \in \R_+^*$ and one has: $\div(f)=D$.
\end{prop}
\proof By applying the equation \eqref{propthet}, one has 
$$
f(p\lambda)-f(\lambda) =\sum_i(h_i\lambda -\mu_i)-\sum_j(h'_j\lambda -\mu'_j)-(p-1)h\lambda=0
$$
since by hypothesis $\sum_i \mu_i-\sum_j\mu'_j=\deg(D)=0$ and $(p-1)h=\sum h_i-\sum h'_j$. Moreover, the divisor of $f$ is equal to $D$ by construction since each theta function $\Theta_{h_i,\mu_i}$ (resp. $\Theta_{h'_j,\mu'_j}$) contributes
with the term $\delta(h_i,\mu_i)$ (resp. $\delta(h'_j,\mu'_j)$) while $-h\lambda $ does not contribute at all.\endproof

\begin{thm}\label{thmtheta1}
Let $f\in \cK(C_p)$. Then $f$ admits the following canonical decomposition 
\begin{equation}\label{thmthet2}
	f(\lambda)=\sum_i \Theta_{h_i,\mu_i}(\lambda)-\sum_j \Theta_{h'_j,\mu'_j}(\lambda)-h\lambda +c
\end{equation}
where $c\in \R$, $(p-1)h=\sum h_i-\sum h'_j$ and $h_i\leq \mu_i<ph_i$, $h'_j\leq \mu_j<ph'_j$.	
\end{thm}
\proof Let $D=(f)$ be the principal divisor of $f$, and $D=D_+-D_-$ its decomposition with 
$(h_i,\mu_i)\in H_p^+\times \R_+^*$, $(h'_j,\mu'_j)\in H_p^+\times \R_+^*$ such that $D_+=\sum \delta(h_i,\mu_i)$ and $D_-=\sum \delta(h'_j,\mu'_j)$. Since $\delta(ph,\mu)=\delta(h,\mu)$ we can choose the $h_i$ and $h'_j$ in such a way that $\mu h^{-1}\in [1,p)$ for all of them, which gives $h_i\leq \mu_i<ph_i$, $h'_j\leq \mu_j<ph'_j$. Proposition \ref{proptheta1} then shows that the function defined in \eqref{propthet2} belongs to $\cK(C_p)$ and has the same divisor as $f$ thus it differs from $f$ by a constant $c\in\R$ and one gets \eqref{thmthet2}. \endproof 
\subsection{Riemann-Roch theorem of type II}

In this subsection we define a Riemann-Roch problem for divisors on the curves $C_p$ and prove a Riemann-Roch formula. To this end, we introduce the notion of continuous dimension for the $\mathbb R_{\rm max}$-modules $H^0(C_p,\cO(D))$ associated to divisors $D$ on $C_p$.

It follows from Proposition \ref{cartierdiv} that
the notion of divisors on $C_p$ coincides with the notion of global section of the 
sheaf $\cdiv(C_p)=\cK_p^*/\cO_p^*$ of Cartier divisors. 
In analogy with the classical case,  a Cartier divisor $D$ described by local sections $f_i\in \cK_p^*(W_i)/\cO_p^*(W_i)$, defines a subsheaf $\cO(D)$ of $\cK_p$. This is the sheaf of $\cO_p$-modules generated on $W_i$ by the (multiplicative) inverses of the $f_i$, \ie in the additive notation of this tropical set-up by the functions $-f_i$. The sections of $\cO(D)$ are given locally by the rational functions $f\in \cK_p(W)$ which satisfy the inequality $D+(f)\geq 0$, where $(f)$ denotes the principal divisor associated to $f$. By construction, $\cO(D)$ is a  sheaf of $\cO_p$-modules, and in particular its global sections define a module over $\rma$
\begin{equation}\label{rrproblem}
H^0(D):=\Gamma(C_p,\cO(D)) =\{f\in \cK(C_p)\mid D+(f)\geq 0\}.
\end{equation}
It follows from \eqref{ordcomp1} that $f,g\in H^0(D)\Rightarrow f\vee g\in H^0(D)$. The constant function $-\infty$ 
 is, by convention, contained in all the $H^0(D)$ and plays the role of the $0$-element. We use the notation $H^0(D)=0$, rather than $H^0(D)=\{-\infty\}$, to mean that  $H^0(D)$ does not contain any other  $f$.

\begin{lem}\label{periodrr} $(i)$~If $\deg(D)<0$ then $H^0(D)=0$.  \newline
$(ii)$~If $\deg(D)>0$ then $H^0(D)\neq 0$.
\end{lem}
\proof $(i)$~The condition $D+(f)\geq 0$ implies that $\deg(D+(f))\geq 0$ and hence $\deg(D)\geq 0$.

$(ii)$~Assume $\deg(D)=\lambda >0$, set $H=\lambda H_p$ and let $P_\lambda$ be the positive divisor which vanishes for $H'\in C_p$,  $H'\neq H$ and takes the value $\lambda$ at $H$. By construction $\lambda\in H$  and $P_\lambda$ is an effective (\ie positive) divisor such that $\deg(P_\lambda)=\lambda$. Moreover one has $\chi(P_\lambda)=1\in \Z/(p-1)\Z$. Let $m\in \{1,\ldots, p-1\}$ be an integer congruent to $\chi(D)$ mod. $p-1$, then it follows that the divisor $D'=D-mP_{\lambda/m}$ fulfills $\deg(D')=0$ and $\chi(D')=0$ and  hence it is principal. Let $f\in \cK_p$ with $D'+(f)=0$. One has $D+(f)=mP_{\lambda/m}\geq 0$ and thus $f\in H^0(D)\neq 0$.\endproof

The slopes of the functions $f\in \cK(C_p)$ can have arbitrarily small real size since the group $H_p\subset \R$ is dense. This shows that when $\deg(D)>0$, \eqref{rrproblem} will in general yield an infinite dimensional space of solutions. However, notice that the group $H_p$ is discrete when embedded diagonally in $\Q_p\times \R$, and this fact allows one to obtain a natural norm on  sections of $\cK_p$ by implementing the $p$-adic norm of the slopes. 
In the following, we choose to normalize the $p$-adic norm $\vert h\vert_p\geq 0$, $\forall h\in H_p$   so that $\vert p\vert_p=1/p$. 

Let $f$ be a continuous, piecewise affine function on $\R_+^*$ with slopes $h_\pm(u)\in H_p$   and such that $f(pu)=f(u)$. The slope $h_\pm(p\lambda)$ of $f$ at $p\lambda$ is $h_\pm(\lambda)/p$ since
$$
h_\pm(p\lambda):=\lim_{\epsilon\to 0\pm} \frac{f(p\lambda+\epsilon)-f(p\lambda)}{\epsilon}=
\lim_{\delta\to 0\pm} \frac{f(p\lambda+p\delta)-f(p\lambda)}{p\delta}=\frac 1p 
\lim_{\delta\to 0\pm}  \frac{f(\lambda+\delta)-f(\lambda)}{\delta}.
$$
The value of $\vert h_\pm(\lambda)\vert_p/\lambda$ is unchanged if one replaces $\lambda$ by $p\lambda$ since  the $p$-adic norm of $h_\pm(\lambda)/p$ is $p\vert h_\pm(\lambda)\vert_p$.
\begin{defn}\label{pnorm} Let $f\in \cK(C_p)$. We set  	
\begin{equation}\label{defnnp}
\Vert f\Vert_p :=\max \{\vert h(\lambda)\vert_p/\lambda\mid  \lambda\in \R_+^*\}\end{equation}
where $h(\lambda)\in H_p$ is the\footnote{at a point of discontinuity of the slopes one takes the max of the two values $\vert h_\pm(\lambda)\vert_p/\lambda$ in \eqref{defnnp}} slope of $f$ at $\lambda$, and $\vert h(\lambda)\vert_p$ its $p$-adic norm.
\end{defn}
  One has the following compatibility with the semiring structure of $\cK(C_p)$
\begin{prop}\label{pnormcomp} Let  $f,g\in \cK(C_p)$. One has\newline
$(i)$~$\Vert f\vee g\Vert_p\leq \max\{\Vert f\Vert_p,\Vert g\Vert_p\}$.  \newline
$(ii)$~$\Vert f + g\Vert_p\leq \max\{\Vert f\Vert_p,\Vert g\Vert_p\}$.\newline
$(iii)$~$\Vert p^a f\Vert_p=p^{-a}\Vert f\Vert_p$. \newline
$(iv)$~$\Vert f\Vert_p\leq 1$ iff the restriction of $f$ to $[1,p]$ has integral slopes. \newline
$(v)$~Let $D\in \div(C_p)$ be a divisor. The following formula defines an increasing filtration on $H^0(D)$ by $\rma$-submodules
\begin{equation*}\label{defnfilt}
H^0(D)^\rho:=\{f\in H^0(D)\mid \Vert f\Vert_p\leq \rho\}.
\end{equation*}
\end{prop}
\proof $(i)$~At a point $\lambda\in \R_+^*$ the set of slopes (there can be two) of $f\vee g$ is a subset of the union of the sets of slopes of $f$ and $g$, thus one obtains the stated inequality.

$(ii)$~At a point $\lambda\in \R_+^*$, the slope of $f+g$ is the sum of the slopes of $f$ and $g$ and the ultrametric inequality for the $p$-adic norm gives the required result.

$(iii)$~The equality follows from the equality $\vert p^a x\vert_p=p^{-a}\vert x\vert_p$.

$(iv)$~Using the invariance of $\vert h(\lambda)\vert_p/\lambda$ under $\lambda\mapsto p\lambda$, one has
$$
\Vert f\Vert_p :=\max \{\vert h(\lambda)\vert_p/\lambda\mid \lambda \in [1,p]\}.
$$
Note that the value $\vert h(1)_-\vert_p$ corresponding to the ingoing slope at $\lambda=1$ is taken into account for $\lambda\in[1,p]$ near $p$ as the limit of the values $\vert h(\lambda)\vert_p/\lambda$ when $\lambda\to p$, since $h(p)_-=h(1)_-/p$.
If the restriction of $f$ to $[1,p]$ has integral slopes, one has $\vert h(\lambda)\vert_p\leq 1$ for all $\lambda\in[1,p]$ and thus $\Vert f\Vert_p\leq 1$. Conversely, if $\Vert f\Vert_p\leq 1$ one has $\vert h(\lambda)\vert_p<p$ for all $\lambda\in[1,p)$ and thus, since $h(\lambda)\in H_p$, one gets $h(\lambda)\in \Z$ for all $\lambda\in[1,p)$ as required.

$(v)$~It follows from $(i)$ and $(ii)$ that $H^0(D)^\rho$ is an $\rma$-submodule of $H^0(D)$, moreover one easily sees that  $H^0(D)^\rho\subset H^0(D)^{\rho'}$ for $\rho<\rho'$.
 \endproof
 
 The next step is to define the continuous dimension of the module $H^0(D)$ using the filtration by the submodules $H^0(D)^\rho$ by means of a formula of the form 
$$
\cdim(H^0(D)):=\lim_{\rho\to \infty}\frac 1 \rho \dim(H^0(D)^\rho),
$$
where, on the right hand side,  one uses a suitable notion of integer valued dimension for $\rma$-modules. In our context, the most natural notion is 
 that of  ``topological dimension" that counts the number of real parameters on which a general element depends.  The original  definition for such dimension is due to Lebesgue. We use the reference \cite{Pears}.

\begin{defn}\label{defntopdim} Let $X$ be a topological space. The {\em topological dimension} $\tdim(X)$ of $X$ is the smallest integer $n$ such that every open cover $\cU$ of $X$ admits a refinement $\cV$ such that every point of $X$ is in at most $n+1$ elements of $\cV$.
\end{defn}
When working with the modules $H^0(D)$ we use the topology of uniform convergence for $\rma$-valued functions on the interval $[1,p]$ (or equivalently by periodicity on $\R_+^*$). The distance defining this topology is given by 
\begin{equation}\label{distop}
	d(f,g)= \max_{x\in [1,p]} \vert f(x)-g(x)\vert.
\end{equation}
We define
\begin{equation}\label{rr1}
\cdim(H^0(D)):=\lim_{n\to \infty} p^{-n}\tdim(H^0(D)^{p^n})
\end{equation}
Our next goal is to prove that the above limit exists and one obtains a Riemann-Roch formula.

\begin{thm}\label{RRperiodic}
$(i)$~Let $D\in \div(C_p)$ be a divisor with $\deg(D)\geq 0$. Then the limit in \eqref{rr1} converges and one has  
$\cdim(H^0(D))=\deg(D)$.\newline
$(ii)$~The Riemann-Roch formula holds
\begin{equation*}\label{rr2}
\cdim(H^0(D))-\cdim(H^0(-D))=\deg(D)\qquad \forall D\in \div(C_p).
\end{equation*}	
\end{thm}
The proof of Theorem \ref{RRperiodic} will be given later and follows by combining  Lemma \ref{lemdivrel} below with the  understanding of $H^0(D)$ in the special case  discussed in Lemma \ref{Enp}.

\begin{lem}\label{lemdivrel} Let $D\in \div(C_p)$ and $f\in \cK(C_p)$.  Then\newline
$(i)$~For $D'=D+(f)$ and  for any non-negative integer $n$ such that $\Vert f\Vert_p\leq p^n$, the map $H^0(D)\to H^0(D')$,  $\xi\mapsto \xi-f$ induces an isomorphism of $H^0(D)^{p^n}$ with $H^0(D')^{p^n}$.\newline
$(ii)$~The absolute Frobenius $f\mapsto p^n f$ determines a twisted isomorphism  of $\rma$-modules
$$
F_{p^n}:H^0(D)^{p^n}\to H^0( p^n D)^{1}.
$$ 
This map preserves the topological dimension. 
\end{lem}
\proof $(i)$~By Lemma \ref{pnormcomp}, $(ii)$, one has $\Vert \xi-f\Vert_p\leq \max\{\Vert \xi\Vert_p,\Vert f\Vert_p\}$ so that if $\Vert f\Vert_p\leq p^n$ one derives
$$
\xi \in H^0(D)^{p^n}\iff \Vert \xi\Vert_p\leq p^n\iff \Vert \xi-f\Vert_p\leq p^n\iff
\xi -f\in H^0(D')^{p^n}.
$$
$(ii)$~The twisting is by the Frobenius $\fr_{p^n}\in \Aut(\rma)$ and occurs since $F_{p^n}(f+a)=F_{p^n}(f)+p^n a$ for any scalar $a\in \rma$. This operation does not affect the topological dimension. The principal divisor $(p^nf)$ is equal to $p^n\times (f)$ and one has, using Proposition \ref{pnormcomp} $(iii)$ 
$$
D+(f)\geq 0\iff p^n D+(p^nf)\geq 0, \ \ \Vert f\Vert_p\leq p^n\iff \Vert p^n f\Vert_p\leq 1.
$$
Thus the map $f\mapsto p^n f$ gives a twisted isomorphism of $\rma$-modules as stated. 
\endproof 

Next, we determine the topological dimension of the $\rma$-module $\cE_{N,p}:=H^0(D)^{1}$ which is associated to the divisor $D:=(H_p,N)$, where $N>0$ is an  integer.

\begin{lem}\label{Enp} $(i)$~The module $\cE_{N,p}$ is the $\rma$-module of 
convex (continuous), piecewise affine functions on $[1,p]$ with integral slopes, such that  $f(1)=f(p)$ and $-f'_+(1)+pf'_-(p)\leq N$.\newline
$(ii)$~Let $a\in \{1,\ldots,N-p\}$ and denote by $b=E((N-a)/p)\geq 1$ the integer part of $(N-a)/p$.  Let 
$\phi_a(x):=\max\{-a(x-1),b(x-p)\}$, for $x\in [1,p]$. Then $\phi_a\in \cE_{N,p}$.\newline
$(iii)$~For $\epsilon > 0$, let $\Delta_{N-p}^\epsilon:=\{(t_1,\ldots,t_{N-p})\mid 0<t_1<\ldots <t_{N-p}<\epsilon\}$. The following map  is continuous and injective  for $\epsilon$ sufficiently small ($\phi_0:=0$ by convention)
  $$h:\R\times \Delta_{N-p}^\epsilon\to \cE_{N,p},\qquad h(t_0,\ldots ,t_{N-p}):=\vee_0^{N-p} (\phi_{N-p-j}-\sum_0^jt_i).$$ 
  $(iv)$~$\tdim(\cE_{N,p})=N-p+1$.
\end{lem}
\proof 
 $(i)$~By Proposition \ref{pnormcomp} $(iv)$, the condition $\Vert f\Vert_p\leq 1$ means that  the restriction of $f$ to $[1,p]$ has integral slopes. The condition $D+(f)\geq 0$ means that $f$ is convex, piecewise affine, (continuous) inside $[1,p]$ and that $N+\ord_{H_p} f\geq 0$. These properties imply $(i)$ since
$$
\ord_{H_p} f=f'_+(1)-pf'_-(p), \qquad N+\ord_{H_p} f\geq 0\iff  -f'_+(1)+pf'_-(p)\leq N.
$$
$(ii)$~By construction, the function $f=\phi_a$ is continuous, convex, piecewise affine with integral slopes and $\phi_a(1)=\phi_a(p)=0$. Moreover, at the point $x=\frac{a+b p}{a+b}$ where $-a(x-1)=b(x-p)$, the slope of $\phi_a$ changes from $-a$ to $b$. The point $x$ is inside the interval $(1,p)$ since $a>0,b>0$. Thus the slope of $\phi_a$  is $-a$ near $1$ and $b$ near $p$. The condition $-f'_+(1)+pf'_-(p)\leq N$ is fulfilled since $a+pb\leq N$. 
\begin{figure}
\begin{center}
\includegraphics[scale=0.7]{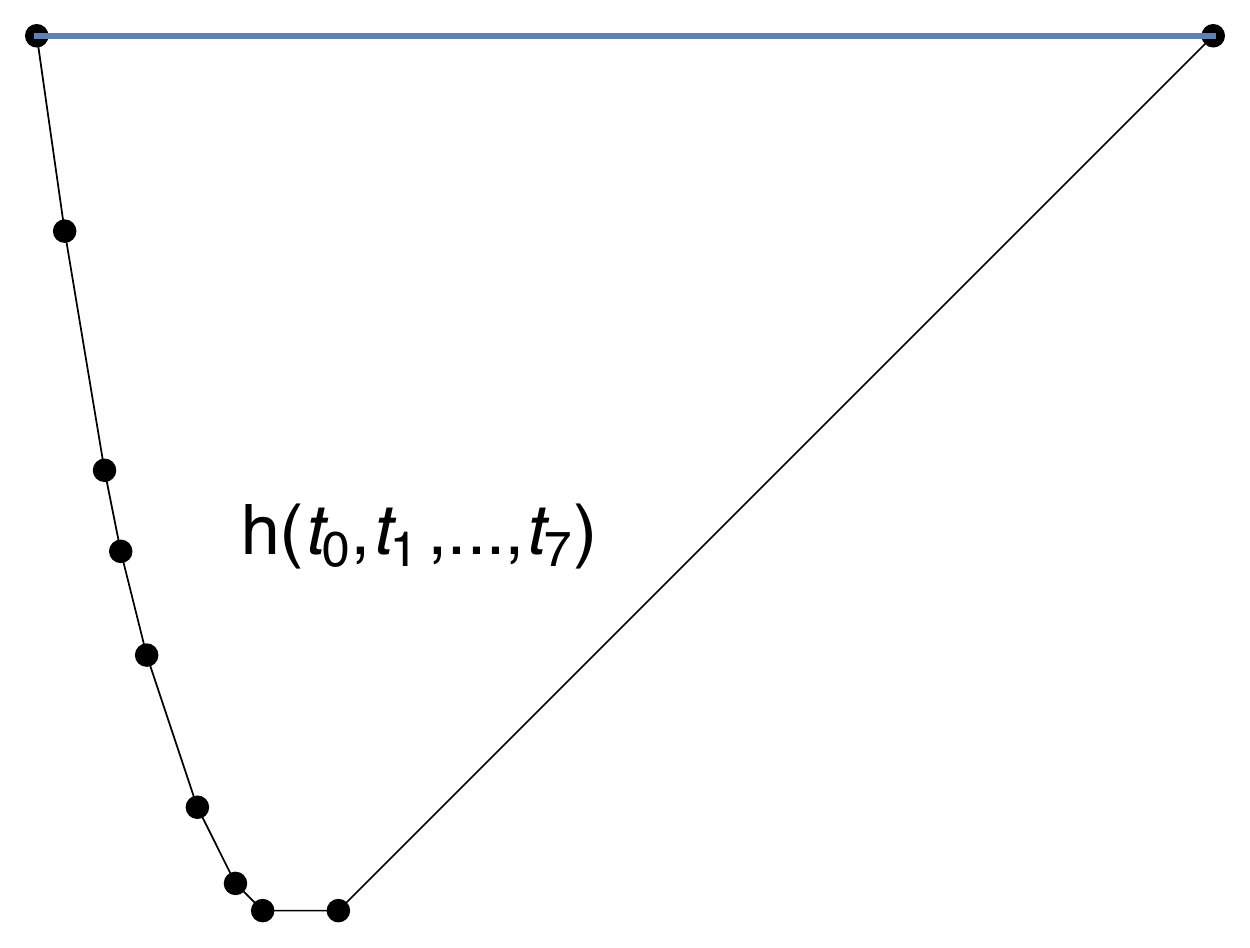}
\end{center}
\caption{The function $h(t_0,\ldots ,t_{N-p})$ with $N-p=7$. \label{trop3} }
\end{figure}
\newline
$(iii)$~Let $\epsilon \leq \frac{p-1}{N-p+1}$. Then on $[1,1+\epsilon]$ all the $\phi_a(x)$'s coincide with $-a(x-1)$ since with $b=E((N-a)/p)\geq 1$ and $a\in \{1,\ldots,N-p\}$ one has  
$$
\frac{a+b p}{a+b}-1=\frac{b (p-1)}{a+b}=(p-1)(1+a/b)^{-1}\geq \frac{p-1}{N-p+1}.
$$
For $0\leq j<N-p$, let $\psi_j:=\phi_{N-p-j}-\sum_0^jt_i$ and let $\psi_{N-p}=-\sum_0^{N-p}t_i$ be a constant function. 
For $0\leq j\leq N-p$ one has
 $$\psi_j(x)=(j+p-N)(x-1)-\sum_0^jt_i \qqq x\in [1,1+\epsilon].$$ For $0\leq j<N-p$, the $x$ coordinate of the point where the line $$L_j:=\{(x,y)\mid y=(j+p-N)(x-1)-\sum_0^jt_i\}$$ meets $L_{j+1}$ is $1+ t_{j+1}\in [1,1+\epsilon]$.  Thus one has
$$
\psi_j(x)\geq \psi_{j+1}(x)\qqq x\in [1,1+t_{j+1}], \ \ 
\psi_j(x)\leq \psi_{j+1}(x)\qqq x\in [1+t_{j+1},1+\epsilon].
$$
Thus $h(t_0,\ldots ,t_{N-p})(x)=\psi_0(x)\qqq x\in [1,1+t_{1}]$ and for $1\leq j\leq N-p-1$,  $h(t_0,\ldots ,t_{N-p})(x)=\psi_j(x)$, for $x\in [1+t_j,1+t_{j+1}]$ (\cf~Figure~\ref{trop3}). It follows that  $h(x)$ passes from the slope $j+p-N$ to the slope $j+1+p-N$ at the point $1+ t_{j+1}\in [1,1+\epsilon]$. This statement still holds for $j=0$ and  $j=N-p-1$. In the latter case, one needs to check that in a small interval $[1+t_{N-p},1+t_{N-p}+\delta]$, $h=\psi_{N-p}$  and this follows as $t_{N-p}<\epsilon$.  The value of the function $h(t_0,\ldots ,t_{N-p})$ at $\lambda=1$ is $-t_0$, thus one recovers all the parameters $t_j$, $0\leq j\leq N-p$ from the function $h(t_0,\ldots ,t_{N-p})$ and this shows that the map $h:\R\times \Delta_{N-p}^\epsilon\to \cE_{N,p}$ is injective.  \newline 
 $(iv)$~Let $k=\tdim(\cE_{N,p})$. The above construction of $h(t_0,\ldots ,t_{N-p})$ shows that $k\geq N-p+1$. Let $f\in\cE_{N,p}$ be non-constant,  then one has $\beta=f_-'(p)\geq 1$ and $\alpha=-f_+'(1)>0$ while $\alpha+p\beta\leq N$. To determine   the topological dimension of $\cE_{N,p}$ one can fix the values of $\alpha$ and $\beta$ and count on how many real parameters the element $f\in \cE_{N,p}$ depends. The possible values of the slope of $f$ in the interval $[1,p]$ are in the interval 
$[f'(1),f'(p)]= [-\alpha,\beta]$ and $f$ is determined, up to an additive constant, by the points where the slope changes, which gives at most $\alpha+\beta$ points. Since the last of these points is determined by the previous ones in view of the periodicity of $f$, one sees that the number of free  parameters for the turning points is $\alpha+\beta-1$. Since this argument determines $f$ up to an additive constant, $f$ depends on at most $\alpha+\beta$ real parameters. Moreover $\alpha+\beta= \alpha+p\beta-(p-1)\beta\leq N-(p-1)$. Thus the topological dimension $k$ of $\cE_{N,p}$ is less than $\alpha+\beta\leq N-p+1$ and since $k\geq N-p+1$ one gets the equality. 

In Appendix \ref{appA} we shall provide a more detailed description of the $\rma$-module $\cE_{N,p}$ making the above qualitative argument more precise. \endproof 

\proof {\it (of Theorem \ref{RRperiodic})}\newline
$(i)$~Assume first that $\deg(D)=0$. Then one has $H^0(D)=0$ except when there exists a non-trivial solution to $D+(f)\geq 0$. In that case $D$ is equivalent to $0$ and $H^0(D)$ consists of the constant functions which form the module $\rma$ whose topological dimension is $1$. Thus 
the formula \eqref{rr1} gives in all cases that $\cdim(H^0(D))=0$. Assume now that $\delta=\deg(D)>0$ and let $\epsilon>0$. We show that 
$$
\underline{\lim}_{n\to \infty}p^{-n}\tdim(H^0(D)^{p^n})\geq \delta-\epsilon, \ 
\overline{\lim}_{n\to \infty}p^{-n}\tdim(H^0(D)^{p^n})\leq \delta+\epsilon
$$
where $\underline{\lim}$ and $\overline{\lim}$ denote respectively the lim inf and the lim sup. Let $\alpha_j\in H_p$ ($j=1,2$), $\alpha_j>0$, be such that $\delta-\epsilon\leq \alpha_1< \delta$,  $\delta< \alpha_2\leq \delta+\epsilon$ and that $\chi(\alpha_j)=\chi(D)$. Then let (as in Lemma \ref{periodrr}, $(ii)$)  $P_j\in \div^+(C_p)$ be positive divisors such that $\chi(P_j)=0$ with $\deg(P_1)=\delta-\alpha_1$, $\deg(P_2)=\alpha_2-\delta$. One has $\deg(\alpha_1\{1\}+P_1)=\delta$ and $\chi(\alpha_1\{1\}+P_1)=\chi(D)$. Thus     by applying Theorem  \ref{thmjaccp} we get the existence of a function  $f_1\in \cK(C_p)$ such that $\alpha_1\{1\}+P_1=D+(f_1)$ and a function $f_2\in \cK(C_p)$ such that $D+(f_2)+P_2= \alpha_2\{1\}$. Using the natural injective maps, isometric for the distance \eqref{distop}, given by the inclusions 
$$H^0(\alpha_1\{1\})\subset H^0(\alpha_1\{1\}+P_1)=H^0(D+(f_1))\stackrel{+f_1}{\to}H^0(D)$$
one obtains, by Lemma \ref{lemdivrel} $(i)$
$$\tdim(H^0(\alpha_1\{1\})^{p^n})\leq \tdim(H^0(D)^{p^n})$$
as soon as $\Vert f_1\Vert_p\leq p^n$.
 Similarly, one derives
$$
\tdim(H^0(D)^{p^n})=\tdim(H^0(D+(f_2))^{p^n})\leq \tdim(H^0(\alpha_2\{1\})^{p^n}) 
$$
as soon as $\Vert f_2\Vert_p\leq p^n$.
Thus the result will follow provided we show that for any $\alpha\in H_p$, $\alpha>0$ one has 
$$
\lim_{n\to \infty} p^{-n}\tdim(H^0(\alpha\{1\})^{p^n})=\alpha.
$$
By Lemma \ref{lemdivrel} $(ii)$, one has 
$\tdim(H^0(\alpha\{1\})^{p^n})=\tdim(H^0(\alpha p^n\{1\})^{1})$ and for $n$ large enough so that $\alpha p^n$ is an integer, Lemma \ref{Enp} $(iv)$ shows that $\tdim(H^0(\alpha p^n\{1\})^{1})=\alpha p^n -p+1$. Thus one gets as required 
$$
\lim_{n\to \infty} p^{-n}\tdim(H^0(\alpha\{1\})^{p^n})=\lim_{n\to \infty} p^{-n} (\alpha p^n -p+1)=\alpha.
$$
$(ii)$~Follows from $(i)$ and the fact that $H^0(D)=0$ if $\deg(D)<0$. \endproof

\appendix
\appendixpage\addappheadtotoc
%\addappheadtotoc

\section{The structure of the $\rma$-module $\cE_{N,p}$} \label{appA}
In this appendix we provide further properties of the $\rma$-module $\cE_{N,p}$ of convex, piecewise affine, continuous functions $f$ on $[1,p]$ with integral slopes  such that $f(1)=f(p)$ and $-f'(1)+pf'(p)\leq N$. 
The convexity of $f$ implies that its graph is above its tangent at any point \ie one has the inequality
\begin{equation}\label{convexin0}
 f(u)\geq f(v)+f'(v)(u-v)\qqq u,v \in [1,p].
\end{equation}
If $v$ is a point of discontinuity of the derivative, one can take for $f'(v)$ any value between the left and right limits $f'_-(v)\leq f'_+(v)$.

We recall the definition of the extremal rays of an $\rma$-module.
\begin{defn}\label{defnextr} Let $\cE$ be an $\rma$-module, and $f\in \cE$. Then $f$ is {\em extremal} iff the equality $f_1\vee f_2=f$ implies that one of the $f_j=f$.  The $\rma$-submodule $\{f+x\mid x\in \rma\}$ is called an extremal ray.
\end{defn}
\begin{figure}
\begin{center}
\includegraphics[scale=0.5]{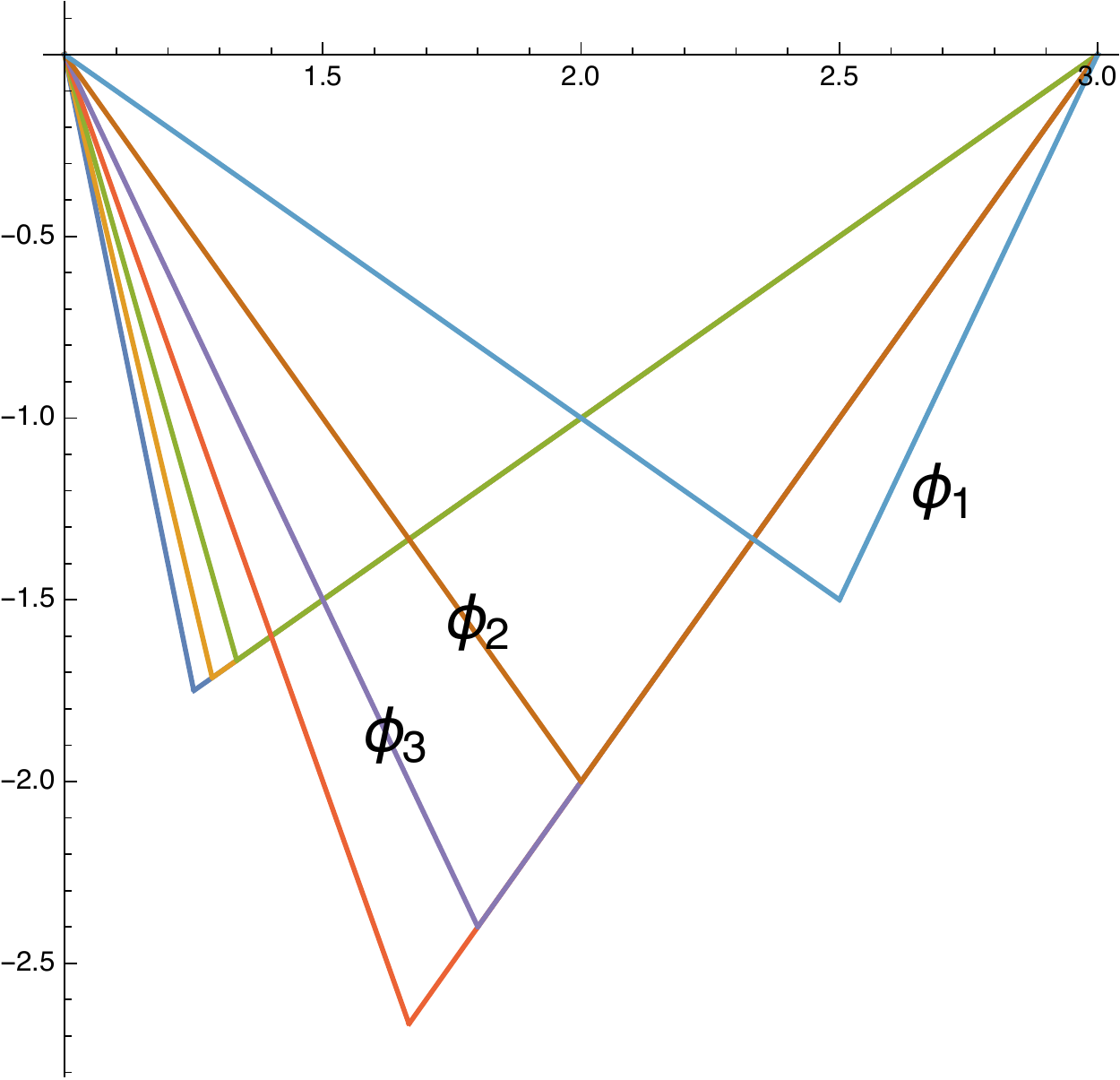}
\end{center}
\caption{The functions $\phi_a$ for $p=3$ and $N=10$; $\phi_0=0$ is constant. \label{trop2} }
\end{figure}
The following proposition describes the structure of the $\rma$-module $\cE_{N,p}$. We use the notation of Lemma~\ref{Enp} (ii).
\begin{prop}\label{pmodgen} The functions  $\phi_a$, for $0\leq a\leq N-p$, are extremal elements of $\cE_{N,p}$ and the elements of the set $\{\phi_a+x\mid x\in \rma\}$ are the extremal rays of $\cE_{N,p}$. These elements generate the module $\cE_{N,p}$, \ie the following map is surjective
\begin{equation}\label{mapsig}
	\sigma:\rma^{N-p+1}\to\cE_{N,p}, \  \  x=(x_a)\mapsto \sigma(x):= \vee (\phi_a+x_a).
\end{equation}	
 \end{prop}
\proof  We  use Lemma \ref{pmodgen0} proved below.  We prove first that $\phi_a$ is extremal. Assume that  $\phi_a=f_1\vee f_2$, with $f_j\in \cE_{N,p}$ for $j=1,2$. Then one of the $f_j$, say $f_1$, agrees with $\phi_a$ in a neighborhood  $V\ni \lambda=1$. By \eqref{dualel1} one has $f_1=\vee(\phi_\ell+ t_\ell)$ for some real numbers $t_\ell$. One of the  $\phi_\ell+ t_\ell$ agrees with $f_1$ in a neighborhood of $1$. Comparing the slopes one gets $\ell=a$ and $t_a=0$. Thus one derives that  $f_1\geq \phi_a$ and since $\phi_a=f_1\vee f_2\geq f_1$ one concludes that $f_1= \phi_a$.   This argument shows that $\phi_a$ is extremal. The equality \eqref{dualel1} shows that any extremal element of $\cE_{N,p}$ is of the form $\phi_a +t$ for some $a$ and $t\in \R$ so that the $\phi_a$ generate all the extremal rays.
Finally \eqref{dualel1} shows that $\sigma$ is surjective.\endproof 
For $a\in \{0,\ldots,N-p\}$, let $\gamma_a:\cE_{N,p}\to \rma$ be defined by 
 \begin{equation}\label{dualel}
\gamma_a(f):=-\max_{x\in [1,p]}  (\phi_a- f)(x).
\end{equation}
\begin{lem}\label{pmodgen0} $(i)$~Let $f\in\cE_{N,p}$ and $x\in (1,p)$. Then, there exists $a\in \{0,\ldots,N-p\}$ such that 
\begin{equation}\label{convextoget}
f(z)-f(x)\geq \phi_a(z)-\phi_a(x),\ \ \forall z\in [1,p].
\end{equation}
$(ii)$~Let $f\in\cE_{N,p}$, then with $\gamma_a(f)$ as in \eqref{dualel}
\begin{equation}\label{dualel1}
f=\vee(\phi_a+\gamma_a(f)).
\end{equation}
\end{lem}
  \proof $(i)$~If $f$ is constant \eqref{convextoget} holds for $a=0$, thus we assume that $f$ is non-constant. Let $\alpha=-f'(1)$ and apply \eqref{convexin0} for $v=1$, $u=p$. Using the equality $f(1)=f(p)$ one has $\alpha\geq 0$, and in fact $\alpha>0$ since $f$ is non-constant. Let $\beta=f'(p)$, by applying \eqref{convexin0} for $v=p$, $u=1$ and using $f(1)=f(p)$ one has $\beta\geq 0$, and in fact $\beta>0$ since $f$ is non-constant. Furthermore one has $\alpha +p \beta\leq N$ since $f\in\cE_{N,p}$. By continuity of the functions involved and the finiteness of the set $\{0,\ldots,N-p\}$, one can assume that $f$ is affine at $x$. Let us first assume that $f'(x)\leq 0$. Let $a=-f'(x)$ so $a\geq 0$. Since $f$ is convex its derivative is non-decreasing and one has $f'(x)\geq f'(1)$ so that $a\leq \alpha$. Moreover by applying \eqref{convexin0} for $v=x$, $u=1$ one has
$$
f(1)\geq f(x)+f'(x)(1-x)=f(x)+a(x-1).
$$ 
By applying \eqref{convexin0} for $v=p$, $u=x$ one derives
$$
f(x)\geq f(p)+f'(p)(x-p)=f(1)+\beta (x-p).
$$ 
Thus 
$$
f(x)\geq f(1)+\beta (x-p)\geq f(x)+a(x-1)+\beta (x-p)
$$
and one gets $a(x-1)+\beta (x-p)\leq 0$, \ie $a(x-1)\leq \beta(p-x)$. One has $\alpha+p\beta\leq N$  and $a\leq \alpha$ so that $a+p\beta\leq N$ and $\beta\leq (N-a)/p$. But $\beta=f'(p)\in \N$ is an integer $\leq (N-a)/p$ and thus it is $\leq E( (N-a)/p)=b$. Then $\beta \leq b$ and one has 
$a(x-1)\leq \beta(p-x)\leq b(p-x)$. By definition $\phi_a(x):=\max\{-a(x-1),b(x-p)\}$ and thus, since  $-a(x-1)\geq b(x-p)$, one gets $\phi_a(x)=-a(x-1)$. By applying \eqref{convexin0} for $v=x$, $u=z$ one has
$$
f(z)\geq f(x)+f'(x)(z-x)=f(x)-a(z-x)=$$ $$=f(x)-a(z-1)+a(x-1)=f(x)-a(z-1)-\phi_a(x).
$$
Thus 
  \begin{equation}\label{convexinbis}
f(z)-f(x)\geq -a(z-1)-\phi_a(x)\qqq z\in [1,p].
\end{equation}
By applying \eqref{convexin0} for $v=p$, $u=z$, one has
$$
f(z)\geq f(p)+f'(p)(z-p)=f(1)+\beta (z-p)\qqq z\in [1,p].
$$
Since $\beta \leq b$ one gets $\beta (z-p)\geq b(z-p)$ so one derives
$$
f(z)\geq f(1)+\beta (z-p)\geq f(1)+b(z-p).
$$
By \eqref{convexinbis} for $z=1$, one has $f(1)-f(x)\geq -\phi_a(x)$ and we thus get 
\begin{equation}\label{convexinbis1}
f(z)\geq  f(x)  -\phi_a(x) +b(z-p) \qqq z\in [1,p]
\end{equation}
which combines with \eqref{convexinbis} to yield
\begin{equation*}\label{convexinbis2}
f(z)-f(x)\geq  \max\{-a(z-1),b(z-p)\}  -\phi_a(x) \qqq z\in [1,p].
\end{equation*}
This is the required inequality \eqref{convextoget}.

  Let us now assume that $f'(x)>0$. Since $f$ is convex, its derivative is non-decreasing and one has $f'(x)\leq f'(p)=\beta$. Let $a=N-pf'(x)$. One has $-pf'(x)\geq -p\beta$ and thus $a\geq N-p\beta$. Since $\alpha +p\beta\leq N$, one gets $a\geq \alpha$. Moreover one has $f'(x)>0$ and thus $a\leq N-p$. The integer  $b:=E((N-a)/p)$ involved in the definition of $\phi_a$ is equal to $f'(x)$ since $N-a=pf'(x)$ by construction. Thus $b=f'(x)$. 
By applying \eqref{convexin0} for $v=1$, $u=x$ one has
$$
f(x)\geq f(1)+f'(1)(x-1)=f(1)-\alpha(x-1).
$$
By applying \eqref{convexin0} for $v=x$, $u=p$ one has
$$
f(p)\geq f(x)+f'(x)(p-x)=f(x)+b(p-x).
$$
Thus, since $f(p)=f(1)$ one gets 
$$
f(x)\geq f(1)-\alpha(x-1)\geq f(x)+b(p-x)-\alpha(x-1)
$$
so that $b(p-x)-\alpha(x-1)\leq 0$, and  $b(p-x)\leq \alpha(x-1)\leq a(x-1)$ since $a\geq \alpha$.
It follows that $\phi_a(x):=\max\{-a(x-1),b(x-p)\}=b(x-p)$. By applying \eqref{convexin0} for $v=x$, $u=z$ one has
$$
f(z)-f(x)\geq f'(x)(z-x)=b(z-x)=b(z-p)-b(x-p)=b(z-p)-\phi_a(x).
$$
Thus 
  \begin{equation}\label{convexinter}
f(z)-f(x)\geq b(z-p)-\phi_a(x)\qqq z\in [1,p].
\end{equation}
By applying \eqref{convexin0} for $v=1$, $u=z$ one has, for $z\in [1,p]$
$$
 f(z)\geq f(1)+f'(1)(z-1)=f(1)-\alpha(z-1)
$$
and since $\alpha\leq a$ one gets $-\alpha(z-1)\geq -a(z-1)$, so that $f(z)\geq f(1)-a(z-1)$. By 
\eqref{convexinter} for $z=p$ one has $f(p)-f(x)\geq -\phi_a(x)$ and since $f(1)=f(p)$ one gets 
$f(1)\geq f(x)-\phi_a(x)$. It follows that for $z\in [1,p]$ one has 
$$
f(z)\geq f(1)-a(z-1)\geq  f(x)-\phi_a(x)-a(z-1)
$$
and hence
\begin{equation}\label{convexinter1}
f(z)\geq   f(x)-\phi_a(x)-a(z-1)\qqq z\in [1,p].
\end{equation}
Thus combining \eqref{convexinter} and \eqref{convexinter1}, one gets
\begin{equation*}\label{convexinter2}
f(z)-f(x)\geq  \max\{-a(z-1),b(z-p)\}  -\phi_a(x) \qqq z\in [1,p]
\end{equation*}
which is the required inequality \eqref{convextoget}.\newline
$(ii)$~For $t=\gamma_a(f)$ one has $\phi_a+t \leq f$, and $\gamma_a(f)$ is the largest value of $t$ for which this happens. 
 By $(i)$ (of this Lemma \ref{pmodgen0}) one sees that  
 at any $x\in (1,p)$ there exists  $a\in \{0,\ldots,N-p\}$ such that \eqref{convextoget}  holds and hence $\phi_a+t \leq f$, for $t=f(x)-\phi_a(x)$ so that $\gamma_a(f)=\max \{t\in\R\mid \phi_a+t \leq f\}\geq f(x)-\phi_a(x)$. One also has $\phi_a(x)+\gamma_a(f)\leq f(x)$, thus one gets $\gamma_a(f)=f(x)-\phi_a(x)$
 and $f(x)=\phi_a(x)+\gamma_a(f)$. It follows that for all $x\in (1,p)$  one has  $f(x)=\vee_{a\in \{0,\ldots,N-p\}}(\phi_a+\gamma_a(f))(x)$.\endproof
 
The section $\gamma$ of the  map $\sigma$ of \eqref{mapsig} has the following interpretation
\begin{equation*}\label{undersgamma}
	\gamma(f)=\max\{g\in \rma^{N-p+1}\mid \sigma(g)=f\}.
\end{equation*}
Indeed, by \eqref{dualel1} one has $\sigma\circ \gamma=\id$, and since $\gamma_a(f)$ is the largest value of $t$ for which $\phi_a+t \leq f$ 
$$
\sigma(x)=f\Rightarrow f= \vee (\phi_a+x_a)\Rightarrow x_a\leq \gamma_a(f)\qqq a
\in \{0,\ldots,N-p\}.
$$
For $f,g\in\cE_{N,p}$ one has $\sigma^{-1}(f\vee g)\supset \sigma^{-1}(f)\vee \sigma^{-1}(g)$. However the equality does not hold in general which explains why $\gamma$ does not preserve $\vee $ in general. Given $f\in\cE_{N,p}$ the fiber $\sigma^{-1}(f)$ is of the form, with $S_f\subset \{0,\ldots,N-p\}$ a suitable subset depending on $f$
$$
\sigma^{-1}(f)=\{x=(x_a)\in \rma^{N-p+1}\mid x_a=\gamma_a(f)\qqq a\in S_f, \ 
x_a\leq\gamma_a(f)\qqq a\notin S_f\}.
$$

\section{Why semirings and characteristic one}\label{appzmax}

It is legitimate to wonder in which sense the algebra of semirings and its specialization to ``characteristic one" (\ie the condition $1+1=1$) is a ``natural" one. We examine this question in this Appendix as an illustration  of the saying ``The devil is in the details".

\subsection{Semiadditive category}

We refer to \cite{MacLane} chapter VIII, exercice 2.4 and recall the following 

\begin{defn}\label{defnsemiadd}  A category $\mathscr C$ is {\em semiadditive} if it has a zero object (both initial and final), all finite products and coproducts exist,  and for any pair of objects $M,N$ of $\mathscr C$ the canonical morphism 
from the coproduct  to the product is an isomorphism
$$
\gamma_{M,N}:M \amalg N\to M\times N, \qquad  \gamma_{M,N}=\left( \begin{array}{cc}
\id_M & 0  \\
0 & \id_N   \end{array} \right).
$$ 
\end{defn}
The object $M \amalg N\simeq  M\times N$ is then called the biproduct of $M$ and $N$, and denoted by  $M\oplus N$. Moreover  
matrix notation works for morphisms. In particular, the matrix $\left( \begin{array}{cc}
\id_M &  \id_M   \end{array} \right)$ gives rise to a canonical element   $s_M\in\Hom_\mathscr C(M\times M,M)$. For any pair of objects $A,M$ of $\mathscr C$, this construction turns $\Hom_\mathscr C(A,M)$ into a commutative monoid with a zero element. Given $\alpha,\beta\in \Hom_\mathscr C(A,M)$, their ``sum" is obtained as the composition
$s_M\circ(\alpha,\beta)$, where the pair $(\alpha,\beta)$ is viewed as an element of $\Hom_\mathscr C(A,M\times M)=\Hom_\mathscr C(A,M)\times \Hom_\mathscr C(A,M)$. Moreover, composition of morphisms is distributive over this addition and thus one has 
\begin{lem}\label{lemsemring} Let $\mathscr C$ be a semiadditive category. Then for any object $M$ of $\mathscr C$ the 
set $\End_\mathscr C(M):=\Hom_\mathscr C(M,M)$ is canonically endowed with the structure of semiring.
\end{lem}

\subsection{The uniqueness of $\B$ and $\zmax$}

Lemma \ref{lemsemring} shows  that semirings occur naturally in semiadditive categories. It is well known that finite fields are commutative rings and their multiplicative group is cyclic. When accepting semirings one obtains two notable complements to the list of finite fields  and both structures are semirings of characteristic one. Indeed, there is only one finite semifield which is not a field, and moreover while there is no field whose multiplicative group is infinite cyclic, there is a unique semifield with this property. More precisely, one has
\begin{prop}	
$(i)$~If a semifield contains a non-trivial root of unity it is a field.\newline
$(ii)$~The only finite semifield which is not a field is $\B$.\newline
$(iii)$~The only  semifield whose multiplicative group is infinite cyclic is $\zmax$.
\end{prop}
\proof $(i)$~By Proposition 4.41 of \cite{Golan}, any element of finite order in a semifield $F$ which is not a field, is equal to $1$. Indeed, with $u^n=1$ one finds that $a:=1+u+\ldots +u^{n-1}$ fulfills $ua=a$ so that $u=1$ if $a\neq 0$. If $a=0$ one gets an element $b:=u+\ldots +u^{n-1}$ such that $1+b=0$ and one has a field. \newline
$(ii)$~By $(i)$~all elements of a finite semifield $F$ which is not a field are $0$ or $1$. Moreover the equation $1+1=1$ holds thus  $F=\B$ is the only finite semifield which is not a field.\newline
$(iii)$~Let $\K$ be a semifield whose multiplicative group is infinite cyclic. The map $\chi:\N\to \K$ defined inductively by $\chi(n+1)=\chi(n)+1$ is a morphism of semirings. Assume first that $\chi$ is injective. Then $\chi(n)\in G$ for all $n\neq 0$, where $G$ is the infinite cyclic multiplicative group of $\K$. Thus there would exist two distinct primes $p,p'$ and integers $a,a'>0$ such that $\chi(p)^a=\chi(p')^{a'}$ but this contradicts the injectivity of $\chi$. Thus $\chi$ is not injective and we let $n<m$ be such that $\chi(n)=\chi(m)$. Let $k\in \N$ and assume $k\geq m$. Write $k=am+b$ with $a>0$, $0\leq b<m$. Then one has $$
\chi(k)=\chi(am+b)=\chi(a)\chi(m)+\chi(b)=\chi(a)\chi(n)+\chi(b)=\chi(an+b)=\chi(k')
$$
where $k'=an+b<am+b=k$. This shows that the image  $\chi(\N)$ is the same as the image
$\chi(\{0,\ldots ,m-1\}$ and is hence finite. Since the multiplicative group of $\K$ is infinite cyclic, it follows that every element of $\chi(\N)$ is in $\{0,1\}\subset \K$. There are two cases. If $1+1=0$, then one is in a field of characteristic $2$ and taking a generator $q$ of the multiplicative group one has an equation of the form $1+q=q^a$ for some $a\notin\{0,1\}$ showing that $\K$ is an algebraic extension of $\F_2$ and thus it is finite. If $1+1=1$ one is in characteristic $1$ and since $\K$ is a semifield the maps 
$\fr_n$, $\fr_n(x):=x^n$, are injective endomorphisms (\cite{Golan}  Proposition 4.43). Let then $q$ be a generator of the multiplicative group and $a\in \Z$ such that $1+q=q^a$. One then has, using the maps $\fr_n$, that 
$1+q^n=q^{na}$ for all $n>0$. For  $a>0$ one has $1+q^a=q^{a^2}$ and the associativity implies 
$$
q^{a^2}=1+q^a=1+(1+q)=(1+1)+q=1+q=q^a
$$
so that $a^2=a$ and $a=1$. For $a =-t\leq 0$ one has $1+q^a=q^a(1+q^t)=q^{a+ta}$ and by associativity one gets as above $1+q^a=1+(1+q)=(1+1)+q=1+q=q^a$ so that $a+ta = a$ and $t=a=0$. Thus the only possible values of $a$ are $0$ and $1$. These values correspond to the two possible choices of the generator of the multiplicative group of $\zmax$ and one obtains the required isomorphism $\K\simeq \zmax$.\endproof

\section{The tropical framework}\label{apptropic}

In this paper we have constructed the scaling site $\scal2$ by implementing the extension of scalars (from $\B$ to $\rmax$) for the arithmetic site $\aarith$. It is rather intriguing   that one would have been led to the same definition of $\scal2$ 
by analyzing the well known results on the localization of zeros of analytic functions which involve Newton polygons in the non-archimedean case and the Jensen's formula in the complex case. Indeed, these results combine to show that the tropical half line, \ie   $(0,\infty)$ endowed with the structure sheaf of convex  piecewise affine functions with integral slopes, gives a  framework, common to all these cases, for the localization of zeros of analytic functions in the punctured unit disk. The additional structure involved in the scaling site, \ie the action of $\nt$ by multiplication on the tropical half-line corresponds, as shown below in \eqref{tropicalizationscal} and \eqref{tropicalizationscalbis}, to the transformation on functions given by the composition with the $n$-th power of the variable.
  The tropical notion of ``zeros" of a convex  piecewise affine function $f$ with integral slope is that a zero of order $k$ occurs at a point of discontinuity of the derivative $f'$ with $k$ equal to the sum of the outgoing slopes. The conceptual meaning of this notion is understood either by using Cartier divisors or by 
implementing  the natural hyperstructure on $\rma$ as shown in \cite{Viro,Viro1}. 
\subsection{Newton polygons}\label{sectnewton}

Let $K$ be a complete and algebraically closed extension of $\Q_p$ and  $v(x):=-\log\vert x \vert$ be the valuation. The tropicalization (see for instance \cite{Eins}) of a series with coefficients in $K$  is obtained by applying the transformation $a\mapsto \log\vert a \vert=-v(a)$ to the coefficients using the simple correspondences
$$
+ \rightarrow \vee , \ \ \times \rightarrow +,\  \ X^n\rightarrow -nx.
$$
In this way a sum $\sum a_n X^n$ of monomials is replaced by the max: $\vee (-nx-v(a_n))$. More precisely
\begin{defn}\label{defntropna}  Let $f(X)=\sum a_n X^n$ be a  Laurent  series with coefficients in $K$ and convergent in an annulus  $A(r_1,r_2)= \{z\in K\mid r_1<\vert z\vert <r_2\}$.  Then, the tropicalization  of $f$ is the real valued function of a real parameter
\begin{equation}\label{tropicalization4}
\tau(f)(x):= \max_n\{-nx -v(a_n) \}\qqq x\in (-\log r_2, -\log r_1).
\end{equation}	
\end{defn}
Up to a trivial change of variables, this notion is well known in $p$-adic analysis where the function $-\tau(-x)$ or rather its graph is called the {\em  valuation polygon} of the series (\cf \cite{Robert} chapter 6, \S 1.6). 
 This polygon is dual to the Newton polygon of the series which is, by definition, the lower part of the convex hull of the points of the real plane with coordinates $(j,v(a_j))$, where $v$ is the valuation. By construction, $\tau(f)(x)$ is finite since using the convergence hypothesis the terms $-nx -v(a_n)$ tend to $-\infty$ when $\vert n\vert\to \infty$. Thus one gets a  convex and piecewise affine function. Moreover one has the multiplicativity property (\cf \cite{Robert} chapter 6, Proposition 1.4.2)
\begin{equation*}\label{tropicalization5np}
\tau(fg)(x)=\tau(f)(x)+\tau(g)(x)\qqq x\in (0, \infty),
\end{equation*} 
and the classical result
 \begin{thm}\label{classical} Let $f(X)=\sum a_n X^n$ be a  Laurent  series with coefficients in $K$, convergent in an annulus  $A(r_1,r_2)= \{z\in K\mid r_1<\vert z\vert <r_2\}$.  	Then the valuations $v(z_i)$ of its zeros $z_i\in A(r_1,r_2)$ (counted with multiplicities) are the zeros (in the tropical sense and counted with multiplicities) of the tropicalization $\tau(f)$ in $(-\log r_2, -\log r_1)$.
\end{thm}
This classical result is, for instance, deduced from the case of Laurent polynomials (see \cite{Dwork}, Lemma 1.6) using the approximation argument of  \cite{Robert}, proof of Theorem 6.2.2.1. 
In particular one can take $r_1=0$, $r_2=1$ so that $A(r_1,r_2)$ is the punctured unit  disk $D(0,1)\setminus \{0\}$. In that case the $\tau(f)$ are convex piecewise affine functions on $(0,\infty)$ and one has the compatibility with the action of $\nt$ on functions by $f(X)\mapsto f(X^n)$ 
\begin{equation}\label{tropicalizationscal}
\tau(f(X^n))(x)=\tau(f)(nx)\qqq x\in (0, \infty),\ n\in \nt.
\end{equation}

\subsection{Jensen's Formula}\label{sectjensen}
For complex numbers, unlike in the non-archimedean case, it is no longer true that for a generic radius $r$, the modulus  $\vert f(z)\vert$ is constant on the sphere of radius $r$ and one replaces \eqref{tropicalization4} by
\begin{defn}\label{defntrop} Let $f(z)$ be a holomorphic function in an annulus  $A(r_1,r_2)= \{z\in \C\mid r_1<\vert z\vert <r_2\}$. Its tropicalization is the function on the interval $(-\log r_2, -\log r_1)$
\begin{equation*}\label{tropicalization4bis}
\tau(f)(x):=\frac{1}{2 \pi}\int_0^{2\pi} \log\vert f(e^{-x+i\theta})\vert d\theta.
\end{equation*}	
\end{defn}
 By construction one still has the multiplicativity
\begin{equation*}\label{tropicalization5}
\tau(fg)(x)=\tau(f)(x)+\tau(g)(x)\qqq x\in (0, \infty).
\end{equation*} 
For $x\in (-\log r_2, -\log r_1)$ such that $f$ has no zero on the circle of radius $e^{-x}$, the derivative of $\tau(f)(x)$ is the opposite of the winding number $n(x)$ of the loop $\theta\mapsto f(e^{-x+i\theta})\in \C^\times$. Thus the function $\tau(f)(x)$ is piecewise affine with integral slopes. When the radius $e^{-x}$ of the circle increases the winding number of the associated loop increases by the number of zeros of $f$ in the intermediate annulus and this shows that the function $\tau(f)(x)$ is convex and fulfills  the Jensen's formula (\cf~\cite{Rudin} Theorem 15.15)  which is seen as the analogue of Theorem \ref{classical}. 
\begin{thm}\label{classical1} Let $f(z)$ be a holomorphic function in an annulus  $A(r_1,r_2)= \{z\in \C\mid r_1<\vert z\vert <r_2\}$ and $z_i\in A(r_1,r_2)$ its zeros counted with their multiplicities. Then the values $-\log\vert z_i\vert$  are the zeros (in the tropical sense and counted with multiplicities) of the tropicalization $\tau(f)$ in $(-\log r_2, -\log r_1)$.
\end{thm}  
In particular one can take $r_1=0$, $r_2=1$ so that $A(r_1,r_2)$ is the punctured unit  disk $D(0,1)\setminus \{0\}$. In that case the $\tau(f)$ are convex piecewise affine functions on $(0,\infty)$ and 
  one has the following compatibility with the action of $\nt$ on functions by $f(z)\mapsto f(z^n)$ 
\begin{equation}\label{tropicalizationscalbis}
\tau(f(z^n))(x)=\tau(f)(nx)\qqq x\in (0, \infty),\ n\in \nt.
\end{equation}
This follows from the equality for periodic functions $h(\theta)$
$$
\frac{1}{2 \pi}\int_0^{2\pi}h(n\theta) d\theta=\frac{1}{2 n\pi}\int_0^{2n\pi}h(\alpha)d\alpha=\frac{1}{2 \pi}\int_0^{2\pi}h(u)du.
$$

\end{document}